\documentclass[11pt,centertags,oneside]{article}
\usepackage{amsmath,amstext,amsthm,amscd,typearea}
\usepackage{amssymb,xcolor}
\usepackage{a4wide}
\usepackage[mathscr]{eucal}
\usepackage{mathrsfs}
\usepackage{typearea}
\usepackage{charter}
\usepackage{pdfsync}
\usepackage{url}
\usepackage{enumitem}
\usepackage[T1]{fontenc}
\usepackage[breaklinks=true,bookmarks=false,pagebackref]{hyperref}

\usepackage{xcolor}

\usepackage[a4paper,width=15.2cm,top=4cm,bottom=4cm]{geometry}

\numberwithin{equation}{section}

\newtheorem{theorem}{Theorem}[section]
\newtheorem{proposition}[theorem]{Proposition}
\newtheorem{corollary}[theorem]{Corollary}
\newtheorem{lemma}[theorem]{Lemma}
\newtheorem{remark}[theorem]{Remark}

\theoremstyle{definition}
\newtheorem{definition}[theorem]{Definition}
\newtheorem{example}[theorem]{Example}

\newcommand{\supp}{{\rm Supp}}

\newcommand{\vol}{\mathop{\mathrm{vol}}}

\newcommand{\ddc}{dd^c}
\newcommand{\dc}{ {d^c}}

\DeclareMathOperator{\ric}{Ric}

\newcommand{\PSH}{{\rm PSH}}

\newcommand{\capK}{\text{cap}}

\newcommand{\C}{\mathbb{C}}

\newcommand{\N}{\mathbb{N}}

\newcommand{\R}{\mathbb{R}}
\renewcommand\P{\mathbb{P}}

\title{\bf K\"ahler--Einstein metrics on quasi-projective manifolds}
\providecommand{\keywords}[1]{\textbf{\textit{Keywords:}} #1}
\providecommand{\subject}[1]{\textbf{\textit{Mathematics Subject Classification 2020:}} #1}

\author{Quang-Tuan Dang $\&$ Duc-Viet Vu}

\newcommand{\Addresses}{{
		\bigskip
		\footnotesize
		\noindent
		\textsc{Quang-Tuan Dang, The Abdus Salam International Centre for Theoretical Physics (ICTP), Str. Costiera, 11, 34151 Trieste, TS, Italy}
		\noindent
		\par\nopagebreak
		\noindent
		\textit{E-mail address}: \texttt{qdang@ictp.it}	\\
		\vskip0.2cm
		\footnotesize
		\noindent
		\textsc{Duc-Viet Vu, University of Cologne, Division of Mathematics, Department of Mathematics and Computer Science, Weyertal 86-90, 50931, K\"oln,  Germany}
		\noindent
		\par\nopagebreak
		\noindent
		\textit{E-mail address}: \texttt{dvu@uni-koeln.de}	
}}
\date{\today}

\begin{document}
	
	\maketitle
	
	\begin{abstract} Let $X$ be a compact K\"ahler manifold and $D$ be a simple normal crossing divisor on $X$ such that $K_X+D$ is big and nef. We first prove that the singular K\"ahler--Einstein metric constructed by Berman--Guenancia is almost-complete on $X \backslash D$ in the sense of Tian--Yau. In our second main result, we establish the weak convergence of conic K\"ahler--Einstein metrics of negative curvature to the above-mentioned metric when $K_X+D$ is merely big, answering partly a recent question posed by Biquard--Guenancia. Potentials of low energy play an important role in our approach.   
	\end{abstract}
	\noindent
	\keywords {Monge--Amp\`ere equations, conic K\"ahler--Einstein metrics, almost-complete metrics, domination principle, analytic singularities}.
	\\
	
	\noindent
	\subject{32U15}, {32Q15}.

	\section{Introduction}
	
	Let $M$ be a $n$-dimensional (non-compact) K\"ahler manifold. Let $\eta$ be a closed positive $(1,1)$-current on $M$. Following \cite{Berman-Guenancia,Berman-Boucksom-Jonsson}, we say that $\eta$ has a well-defined Ricci curvature if the non-pluripolar self-product $\langle \eta^n \rangle$  of $\eta$ is well-defined on $M$  (see~\cite{BEGZ} for this notion) and for every local holomorphic volume form $\Omega$ on $M$, one can write $\langle \eta^n\rangle= e^{-2f} i^{n^2}\Omega \wedge \overline \Omega$ for some function $f \in L^1_{\rm loc}$. In this case, we define $\ric \eta:= \frac{i}{\pi} \partial \bar \partial f$. The current $\eta$ is said to be \emph{singular K\"ahler--Einstein metric} on $M$ if $\ric \eta= \lambda \eta$ for some constant $\lambda \in\R$. 
	We are interested in studying  (singular) K\"ahler--Einstein metrics in $M$. We will focus on the most common case where $M$ is the complement of a divisor in a compact K\"ahler manifold. Let us enter details.  
	
	Let $(X,\omega)$ be a compact K\"ahler manifold of dimension $n$.
	Let $D$ be a normal simple crossing divisor on $X$. By abuse of notation, we also denote by $D$ the support of $D$.  About forty years ago, Cheng--Yau \cite{Cheng-Yau} established some necessary conditions for the existence of complete K\"ahler--Einstein metrics on non-compact manifolds (see also~\cite{Mok-Yau} for refinements). Based on this work, Kobayashi \cite{Kobayashi-KE} proved that there exists a complete K\"ahler--Einstein metric of negative curvature (with Poincar\'e singularities near $D$) on $X \backslash D$ if $K_X +D$ is ample. The uniqueness of such a metric follows from~\cite{Yau-schwarz}. 
	
	More generally, Tian--Yau \cite{Tian-Yau-complete-ssurvey} proved the existence and uniqueness of an almost complete K\"ahler--Einstein metric $\Tilde{\omega}_D$ on $X \backslash D$ provided that $K_X+D$ is big and nef, and $K_X+ D$ is ample modulo $D$, (cf. Definition~\ref{def-almost-complete} below for almost complete metrics). 
	 We recall that $K_X+D$ is ample modulo $D$ if for every effective reduced curve $C$ on $X$ not contained in $D$, $(K_X+D)\cdot C>0$, hence
 the non-K\"ahler locus of the Chern class of $K_X+D$ must lie in the support of $D$. We refer to \cite{Tian-Yau-complete-ssurvey} for more information and applications of this metric and to \cite{bando90einstein,auvray13metrics,auvray17asymptotic,auvray17space,Biquard-Guenancia,DiNezza-Lu-quasiprojective,Guenancia-ampleKX,Guenancia-Wu,WuDamin-quasi,WuYau20invariant,
		Fang-Fu,DatarFuSong23} for related results.  
	
	In another  approach, Berman--Guenancia~\cite{Berman-Guenancia} proved that 
	there exists a unique closed positive current $\widehat\omega_D$ of full Monge--Amp\`ere mass (in the sense of~\cite{BEGZ}) in $c_1(K_X +D)$ such that 
	\begin{equation}\label{eq: BermanGuenancia}
	\ric (\widehat\omega_D)= -\widehat\omega_D+ [D]
	\end{equation}
	provided that $K_X+D$ is big, where $[D]$ denotes the current of integration along $D$. If $K_X+ D$ is additionally nef, then $\widehat\omega_D$ is also smooth on the complement of $D$ in the ample locus of $K_X+D$ (see ~\cite{Berman-Guenancia} or Lemma \ref{le-uniform-Cinfintynorm} below). 
	
	In the case where $K_X+D$ is ample, using the local description of $\tilde{\omega}_D$  and the Monge--Amp\`ere equation satisfied by it (see \cite[Page 407]{Kobayashi-KE}), one sees that $\tilde{\omega}_D$ can be extended trivially through $D$ as a closed positive $(1,1)$-current of full Monge--Amp\`ere mass in $K_X+D$ and satisfies the same equation as $\widehat \omega_D$. Hence $\widehat \omega_D= \tilde{\omega}_D$. 
	In other words, we obtain a global interpretation of the complete metric $\tilde{\omega}_D$ if $K_X+D$ is ample (see also \cite{Guenancia-ampleKX} for more information). 
	Here is our first main result giving an analogue in the case where $K_X+D$ is big and nef. 
	
	\begin{theorem} \label{the-TianYaumetric} Let $D$ be a simple normal crossing divisor such that $K_X+ D$ is big and nef. Let $E$ denote the non-ample locus of $K_X+D$.
		Then there exists a unique almost-complete singular K\"ahler--Einstein metric $\omega_D$ on $X \backslash D$, and this metric satisfies the following properties:
		
		(i)  $\omega_D$ is smooth K\"ahler metric on $X \backslash (D \cup E)$, 
		
		(ii) $\omega_D$ can be extended to a closed positive $(1,1)$-current on $X$ which is of full Monge--Amp\`ere mass in the class $c_1(K_X+D)$ and there holds
		$$\ric \omega_D= - \omega_D+ [D]$$
		as currents on $X$. 
	\end{theorem}
	

	We note that since $\omega_D$ is of full Monge--Amp\`ere mass, one has 
	$$\int_{X \backslash (D \cup E)} \omega_D^n = \int_X \big(c_1(K_X+D)\big)^n,$$
	which is to say that the volume of the singular metric $\omega_D$ on $X \backslash D$ is finite and is equal exactly to the volume of the line bundle of $K_X+D$.  Theorem \ref{the-TianYaumetric} shows that this almost-complete metric on $X \backslash D$ is indeed equal to the above metric $\omega_D$ constructed by Berman--Guenancia, using the comparison principle~\cite[Proposition 4.1]{Berman-Guenancia}. 
	It generalizes the previous result of Tian--Yau when $K_X+D$ is big, nef and $K_X+D$ is ample modulo $D$ (because in this case, one has $E \subset D$). 
	

	We now come to the next part of the introduction in which we discuss a very recent question posed by Biquard--Guenancia  \cite{Biquard-Guenancia} about the degeneration of conic K\"ahler--Einstein metrics. Our result gives a direct construction for the above metric $\omega_D$ as the weak limit of conic K\"ahler--Einstein metrics. 
	
	Let $D$ be a simple normal crossing divisor on $X$ such that $K_X+ D$ is big. 
	By \cite{BEGZ,Kolodziej_Acta,Yau1978}, there exists a unique closed positive current $\omega_\epsilon$ of full Monge--Amp\`ere mass in $c_1(K_X+ (1-\epsilon)D)$ so that  \begin{align}\label{eqKEconic}
	\ric  \omega_\epsilon= - \omega_\epsilon+ (1-\epsilon) [D],
	\end{align}
	for every constant $\epsilon \in (0,1)$ small enough. The equation is understood in the sense of currents. 
	If $K_X+ D$ is additionally ample, then the metric $\omega_\epsilon$ is, in fact, conic near $D$ by \cite{Guenancia-Paun} refining \cite{Brendle, CGP, Jeffres-Mazzeo-Rubinstein} (see also \cite{Guenancia-klt, Li-Sun, Szekelyhidi}). For this reason, we say that $\omega_\epsilon$ is a K\"ahler--Einstein metric with a conic singularity of angle $2\pi\epsilon$ along $D$ even if $K_X+D$ is merely big.

	It is thus natural to study the relation between $\omega_\epsilon$ and  $\omega_D$ as $\epsilon\to 0^+$. Biquard--Guenancia \cite{Biquard-Guenancia} asked whether the metric $\omega_\epsilon$ converges to $\omega_D$ as $\epsilon \to 0$ in an appropriate sense. Under some additional assumptions on the positivity of  $K_X+D$, the local $\mathcal{C}^\infty$ convergence was established in \cite{Biquard-Guenancia,Guenancia-ampleKX,Guenancia-Wu}. Here is our second main result proving the weak convergence of $\omega_\epsilon\to \omega_D$ in the setting where $K_X+D$ is merely big.

	\begin{theorem} \label{the-KE-Gue} 
		Let $X$ be a compact K\"ahler manifold and $D$ be a simple normal crossing divisor such that $K_X+D$ is big. Let $\omega_\epsilon$ be the K\"ahler--Einstein metric solving~\eqref{eqKEconic} for small $\epsilon>0$. Then,
		we have the weak convergence $\omega_\epsilon \to \omega_D$ as $\epsilon \to 0^+$ in the sense of currents. 
	\end{theorem}

	We note that Theorem \ref{the-KE-Gue}  remains true in a slightly more general situation that $D=  \sum_{j=1}^m a_j D_{j}$ is a divisor whose support is simple normal crossing such that $a_j \le 1$ for every $1 \le j \le m$. This can be proved essentially along the same line as in our proof of Theorem \ref{the-KE-Gue}. To simplify the presentation, we only consider a simple normal crossing divisor $D$ as above.


	In our proof, we actually don't need the existence of solutions for~\eqref{eq: BermanGuenancia} from~\cite{Berman-Guenancia} and directly show that the sequence of potentials of $(\omega_\epsilon)_\epsilon$ is convergent in $L^1$ as $\epsilon\to 0^+$ (we indeed prove a much stronger property that the sequence of potentials is decreasing in capacity, see Definition \ref{def: decrease-capa}) and the limit current will solve the equation~\eqref{eq: BermanGuenancia}.  
	
	As far as we know, although it was known when the cohomology class $K_X+ (1-\epsilon)D$ is big and nef (\cite{BEGZ}),  it is still open whether the metric $\omega_\epsilon$ is, in general, a genuine metric on the complement of some proper analytic subset in $X$. Thus, we could not for the moment address the question of local smooth convergence of $\omega_\epsilon \to \omega_D$ as $\epsilon\to 0^+$.   We would like to point out, however, that the weak convergence of a sequence of metrics as currents is usually the first step in proving the $\mathcal{C}^\infty$ convergence locally outside a proper analytic subset in $X$.

	In~\cite{Biquard-Guenancia}, under the assumption that $c_1(K_X+(1-\epsilon)D)$ is semi-positive, the authors were able to use the domination principle to show that a suitably normalized sequence of potentials $(u_\epsilon)_\epsilon$ is essentially decreasing to $u_D$, hence obtaining the desired $L^1$ convergence. In another context where $K_X+D$ is ample, the $L^1$ convergence of potentials of $\omega_\epsilon$ was proved in \cite{Guenancia-ampleKX}, using a variational approach. Both these assumptions are not available in our setting.

	We now comment on our proofs of the main results. The proof of Theorem \ref{the-TianYaumetric} is based on a uniform lower bound for solutions of certain classes of complex Monge--Amp\`ere equations whose proof uses an idea from \cite{DDG-family}. The almost-completeness of $\omega_D$ is deduced from a more or less standard uniform $\mathcal{C}^\infty$-estimates for Monge--Amp\`ere equations. For the uniqueness of $\omega_D$, we need a slightly more general version of Yau's Schwarz lemma for singular K\"ahler--Einstein metric, see Theorem \ref{the-yauschwarz-sing} below for details.

	The proof of Theorem \ref{the-KE-Gue} is much more involved because there is no available uniform  $\mathcal{C}^\infty$-estimates (recall that $K_X+D$ is only big). 
	We will indeed prove a more general result (Theorem \ref{the-stability} below) about the stability of complex Monge--Amp\`ere equations.
	In contrast to some well-known non-collapsing situations where the $L^1$ convergence of potentials of metrics is relatively easy to establish, the difficulty in our present setting (i.e., in Theorem \ref{the-stability}) lies in the fact that the complex Monge--Amp\`ere equation corresponding to (\ref{eqKEconic}) has the measure on the right-hand side of unbounded mass on $X$ if $\epsilon=0$. Because of this reason, the expected $L^1$ convergence was obtained in previous works only under some additional assumption.
	
	The strategy of our proof of Theorem \ref{the-stability} is to show that in the absence of any additional positivity (except the minimal assumption that $K_X+D$ is big), one can still use a sort of domination principle to obtain a quasi-decreasing property of potentials of metrics in consideration, hence, get the desired convergence. The version of the domination principle we need is the quantitative one obtained recently in \cite{Do-Vu_quantitative} (see Theorem \ref{the-stability} in Section~\ref{sect: quantitative}). It is robust enough to allow one to compare quasi-plurisubharmonic functions under relatively relaxed situations.  
	\\

	\noindent
	\textbf{Aknowledgment.} The authors would like to thank Henri Guenancia for fruitful discussions. Duc-Viet Vu is partially supported by the Deutsche Forschungsgemeinschaft (DFG, German Research Foundation)-Projektnummer 500055552 and by the ANR-DFG grant QuaSiDy, grant no ANR-21-CE40-0016. Part of this work was done when Quang-Tuan Dang was visiting the Department of Mathematics and Computer Science at the University of Cologne.
	\\
	
	\noindent\textbf{Conventions.} The notations $\lesssim$, $\gtrsim$ stand for inequalities up to a multiplicative uniform positive constant. We use $C$ for a positive constant, whose value may change from line to line. For a divisor $D$, in some context we simply write $D$
	instead of its support ${\rm Supp}(D)$. Denote ${\rm d}=\partial+\Bar{\partial}$ and $\dc=\frac{i}{2\pi}(\Bar{\partial}-\partial)$ so that $\ddc=\frac{i}{\pi}\partial\Bar{\partial}$.	\\	

	\noindent\textbf{Ethics declarations} The author declares no conflict of interest.

	\section{Quantitative domination principle}
	\label{sect: quantitative}

	In this section, we recall the quantitative domination principle obtained in~\cite{Do-Vu_quantitative}. 
	Let $X$ be a $n$-dimensional compact K\"ahler manifold equipped with a K\"ahler metric $\omega$.  A function $u:X\to\mathbb{R}\cup\{-\infty\}$ is quasi-plurisubharmonic
	(qpsh) if it can be locally written as the sum of a plurisubharmonic function and a smooth function. Let $\theta$ be a smooth closed real $(1,1)$-form.
	We say that $u$ is $\theta$-plurisubharmonic ($\theta$-psh)  if it is qpsh and $\theta_u:=\ddc u+\theta\geq 0$ in the sense of currents.  We let $\PSH(X, \theta)$ denote the set of $\theta$-psh functions on $X$. Recall that
	the cohomology class $\{\theta\}$ is {\em big} if there exists $\rho\in\PSH(X,\theta)$ such that $\ddc\rho+\theta\geq \delta\omega$ for some small constant $\delta>0$.
	
	A $\theta$-psh function $u$ is said to have {\em minimal singularities} if it is less singular than any
	other $\theta$-psh function. 
	Define 
	\[ V_\theta:=\sup\{\varphi\in\PSH(X,\theta): \varphi\leq 0 \}.\]
	One can see that $V_\theta$ is a $\theta$-psh function with minimal singularities.
	Let $u_1,\ldots,u_p$, for $1\leq p\le n$ be $\theta$--psh functions and put $\theta_{u_j}:= \ddc u_j+ \theta$. We recall how to define the non-pluripolar product $\theta_{u_1} \wedge \cdots \wedge \theta_{u_n}$. We write locally $\theta_{u_j}= \ddc v_j$, where $v_j$ is psh. By \cite{BT_fine_87, BEGZ}, one knows that the sequence of positive currents 
	$$\mathbf{1}_{\cap_{j=1}^p \{v_j>-k\}}\ddc \max\{v_1, -k\} \wedge \cdots \wedge \ddc \max\{v_p, -k\}$$ 
	is increasing in $k \in \N$ and converges to a closed positive current, which is independent of the choice of local potentials $v_j$'s. Thus we obtain a well-defined global closed positive current on $X$ which is called the \emph{non-pluripolar product} $\theta_{u_1} \wedge \cdots \wedge \theta_{u_p}$ of $\theta_{u_1}, \ldots, \theta_{u_p}$.  If $p=n$, the resulting positive $(n, n)$-current is a Borel
	measure putting no mass on pluripolar sets. For any $u\in\PSH(X,\theta)$,  the non-pluripolar complex Monge--Amp\`ere measure of $u$ is \[\theta_u^n:= (\ddc u+\theta)^n.\]
	Given a potential $\phi\in\PSH(X,\theta)$, we let $\PSH(X,\theta,\phi)$ denote the set of $\theta$-psh functions $u$ such that $u\leq\phi$. We also define by $\mathcal{E}(X,\theta,\phi)$ the set of $u\in\PSH(X,\theta,\phi)$ of full Monge--Amp\`ere mass with respect to $\phi$, i.e.,
	$\int_X\theta_u^n=\int_X\theta_\phi^n$.
	When $\phi=V_\theta$ we simply denote by $\mathcal{E}(X,\theta)$. 
	The volume of a big class $\{\theta\}$ is defined by
 \[\vol(\{\theta\}):=\int_X\theta_{V_\theta}^n. \]
In this case, a $\theta$-psh function $u$ is said to have {\em full Monge--Amp\`ere mass} if $\int_X\theta_u^n=\vol(\theta)$.
 
 The following is a version of the domination principle.
		\begin{proposition}[Non-quantitative domination principle] \label{prop: standarddomination} Let $u\in\mathcal{E}(X,\theta)$ and $v\in {\rm PSH}(X,\theta)$. If $e^{-v}\theta_v^n\geq e^{-u}\theta_u^n$ then $u\geq v$.
	\end{proposition}
	
	\proof
	For every $a>0$, we set $v_a=\max(u,v-a)$. By assumptions, 
	we have $\theta^n_u\leq e^{-a}\theta^n_v$ on the set $\{u<v-a\}$. 
	The locality of the complex
Monge–-Amp\`ere measure with respect to the plurifine topology and the comparison principle yields
	\[ \int_{\{u<v_a\}}\theta^n_u\leq\int_{\{u<v_a\}}e^{-a}\theta^n_{v_a}\leq \int_{\{u<v_a\}}e^{-a}\theta^n_u \]
	This means that $u\geq v_a$ a.e-$\theta^n_u$.
	By the domination principle~\cite[Proposition 3.11]{Lu-Darvas-DiNezza-mono} we obtain $u\geq v_a$, hence $u\geq v-a$ on $X$. Since the latter holds for every $a>0$, we get $u\geq v$ as desired.
	\endproof
	We now recall a particular case of the quantitative domination principle established in  \cite{Do-Vu_quantitative}. Let $\mathcal{W}^-$ be the set of convex, non-decreasing functions $\chi: \R_{\le 0} \to \R_{ \le 0 }$ such that $\chi(0)=0$ and $\chi(-\infty) = -\infty$. 
	Let $\theta$ be a closed smooth $(1,1)$-form in a big cohomology class.  Let $V_\theta$ be the $\theta$-psh function with minimal singularities defined above. 
	Let $\varrho:= \int_X \theta_{V_\theta}^n>0$, which is the volume of the cohomology class of $\theta$. For $\chi \in \mathcal{W}^-$ and   $u \in \mathcal{E}(X,\theta)$, let 
	$$E^0_{\chi, \theta}(u):= - \varrho^{-1} \int_X \chi(u- V_\theta) \theta_u^n$$
	which is called \emph{the (normalized) $\chi$-energy} of $u$.  
	For every Borel set $E$ in $X$,  recall that the {\em capacity} of $E$ is given by 
	$$\capK(E)=\capK_\omega(E):= \sup\left\{\int_{E}\omega_w^n: w \in \PSH(X,\omega), 0 \le w  \le 1 \right\}.$$

	\begin{theorem}[{\cite[Theorem 4.4]{Do-Vu_quantitative}}]\label{the domination}  (Quantitative domination principle)
		Let $A \ge 1$ be a constant, and let $\theta\leq A\omega$ be a closed smooth real $(1,1)$-form in a big cohomology class, and $\varrho:=\int_X\theta_{V_\theta}^n>0$.  Let $B \ge 1$ be a constant, $u_1, u_2\in \mathcal{E}(X, \theta)$  and  $\tilde{\chi} \in \mathcal{W}^-$ such that $\tilde{\chi}(-1)=-1$ and
		$$E^0_{\tilde{\chi}, \theta}(u_1)+E^0_{\tilde{\chi}, \theta}(u_2) \le B.$$
		Assume that there exists a constant $0\leq c<1$ and a Radon measure $\mu$ on $X$ satisfying
		$$\theta_{u_1}^n\leq c\theta_{u_2}^n+\varrho\mu$$
		on $\{u_1<u_2\}$ and 
		$c_{\mu}:=\int_{\{u_1<u_2\}}d\mu\leq 1$.
		Then there exists a constant $C>0$ depending only on $n, X$ and $\omega$   such that 
		$$\capK_{\omega}\{u_1<u_2-\epsilon\}\leq 
		\dfrac{C\vol(X) (A+B)^2}{\epsilon(1-c)h^{\circ n}(1/c_{\mu})},$$
		for every $0<\epsilon<1$,
		where  $h(s)=(-\tilde{\chi}(-s))^{1/2}$ for every $0 \le s \le \infty$.
		
		In particular, if $c_{\mu}=0$ then $\capK_{\omega}\{u_1<u_2-\epsilon\}=0$ for every $\epsilon>0$
		and $u_1\geq u_2$ on $X$.
	\end{theorem}
	
	The standard domination principle (see, e.g.,~\cite[Proposition 3.11]{Lu-Darvas-DiNezza-mono}) corresponds to the case where $c=0$ and $\mu=0$. 
	We underline that it is crucial for us in applications later that we consider the situation where $\mu \not =0$. 
	
	We continue with some more auxiliary results about the continuity of Monge--Amp\`ere operators. 
	We now work in the local context. Let $\Omega$ be an open subset in $\C^n$ and $d\lambda$ denote the Lebesgue measure on $\mathbb{C}^n$. Let $\PSH(\Omega)$ denote the set of psh functions on $\Omega$.
	
	\begin{definition}\label{def: decrease-capa}
		Let $(w_j)_j$ be a sequence of Borel functions on $\Omega$. 
		We say that $(w_j)_j$ is {\em decreasing in capacity} if for every open subset $U$ in $\Omega$,  every compact $K\Subset U$, and for every constant $\delta>0$, for every $\epsilon>0$ there exists an index $j_\epsilon$ such that  
		$$\capK \big(\{w_j - w_{j'} \le -\delta\} \cap K, U\big) \le \epsilon$$ if $j' \ge j \ge j_\epsilon$.
	\end{definition}
	We remark that the capacity in the above definition is in the sense of Bedford--Taylor~\cite{Bedford_Taylor_82}: for any Borel set $E\subset U$,
 \[\capK(E,U):=\sup\left\{\int_E(\ddc u)^n:u\in\PSH(U),0\leq u\leq 1\right\}. \]

	\begin{proposition}\label{prop: cv_decreasingcapacity} Let $(w_j)_j$ be a sequence of uniformly bounded psh functions such that $(w_j)_j$ is decreasing in capacity, and $w_j$ converges to some bounded psh function $w$ in $L^1_{\rm loc}$ as $j \to \infty$. Then we have  $(\ddc w_j)^k \to (\ddc w)^k$ as $j \to \infty$ for every $1 \le k \le n$.
	\end{proposition}
	
	\proof
	 We argue as in the proof of the standard continuity of Monge--Amp\`ere operators for decreasing sequences; see, e.g.~\cite{Demailly_ag}. Following the same lines in the aforementioned reference, we prove the desired limit by induction and assume now that $(\ddc w_j)^{k-1}$ converges weakly to $(\ddc w)^{k-1}$ as $j \to \infty$. To check the desired assertion for $k$, it suffices to show that $w_j(\ddc w_j)^{k-1}$ converges weakly to $w(\ddc w)^{k-1}$.
	Observe that we already have that any limit current of the sequence $w_j(\ddc w_j)^{k-1}$ is bounded from above by $w(\ddc w)^{k-1}$ by the induction hypothesis and Hartogs' lemma. 
	
	We can, without loss of generality, assume that $-2\leq w_j \le -1$ and $-2 \le w\leq -1$. Denote by $w^\epsilon_j$ the standard regularization of $w_j$ by using convolution; cf.~\cite{Demailly_ag}.
	The problem is local. We can thus assume that $\Omega$ is the unit ball in $\C^n$. Let $\psi(z):= |z|^2-1$. We can assume that all $w_j, w_j^\epsilon$ and $w$ coincide with $A\psi$ outside a compact $K$ in $\Omega$, where $A$ is a  constant big enough.
	We set $\beta:=\ddc\psi>0$, which is the standard K\"ahler form in $\C^n$.  Let $\delta>0$, $\epsilon>0$ be  constants. Choose an index $j_{\epsilon,\delta}$ so that for every $j\geq j_{\epsilon,\delta}$,
	$$\capK(\{w_j-w\le-\delta\} \cap K,\Omega) \le \epsilon$$
	Thus there exists a constant $C>0$ independent of $\delta, \epsilon$ such that for every $j \ge j_{\epsilon,\delta}$ we have
	\begin{align*}
	\int_\Omega w(\ddc w)^{k-1}\wedge \beta^{n-k+1}&\leq \int_\Omega(w_j+\delta)(\ddc w)^{k-1}\wedge\beta^{n-k-1}\\
    &\quad+C\capK(\{w_j-w\le-\delta\} \cap K,\Omega)\\
	&\leq \int_\Omega(w_j^{\epsilon}+\delta)(\ddc w)^{k-1}\wedge\beta^{n-k-1}+C\epsilon\\
	&\leq \int_\Omega w\ddc w_j^{\epsilon}\wedge (\ddc w)^{k-2}\wedge\beta^{n-k+1}+C(\epsilon+ \delta)\\
	& \le \int_\Omega (w_j^{\epsilon}+\delta) \ddc w_j^{\epsilon}\wedge (\ddc w)^{k-2}\wedge\beta^{n-k+1}+C(\epsilon+\delta).
	\end{align*} 
	Note that in the above estimates, we have used integration by parts and the fact that $w$ and $w_j^\epsilon$ both vanish on $\partial\Omega$. Repeating this argument, we obtain
	\begin{align*}
	\int_\Omega w(\ddc w)^{k-1}\wedge \beta^{n-k+1}&\leq \int_\Omega w_j^{\epsilon}(\ddc w_j^{\epsilon})^{k-1}\wedge\beta^{n-k+1}+C(\epsilon+ \delta),
	\end{align*}
	for $j \ge j_{\epsilon,\delta}$ and for some constant $C>0$ big enough but independent of $\epsilon$ and $j$. 
	Letting $j\to \infty$, then  $\epsilon \to 0$, and then $\delta \to 0$ we obtain
	$$ \int_\Omega w(\ddc w)^{k-1}\wedge \beta^{n-k+1} \le  \liminf_{j\to \infty}\int_\Omega w_j(\ddc w_j)^{k-1}\wedge \beta^{n-k+1}.$$
	The proof is thus complete.
	\endproof

	\section{Uniform lower bounds}\label{sect: lowerbound}

	Let $D$ be a simple normal crossing divisor in $X$. Let $h$ be a smooth Hermitian metric on $\mathcal{O}_X(D)$ and $s$ a section of $\mathcal{O}_X(D)$ defining $D$, normalized so that $|s|_h\leq 1$. 
    Let $\alpha$ be a big cohomology class, and let $\theta$ be a smooth representative in $\alpha$. Let $(f_j)_{j\ge 1}$ be an increasing sequence of continuous nonnegative functions converging pointwise to $f \in \mathcal{C}^0(X)$ as $j \to \infty$ such that $f_1 \not \equiv 0$. Let $(\theta_j)_{j\in \mathbb{N}}$ be a sequence of smooth closed $(1,1)$-forms in big cohomology classes converging to $\theta$ in the $\mathcal{C}^0$ topology as $j \to \infty$. Assume that for every $j \ge 0$, there exists $u_j \in \mathcal{E}(X, \theta_j)$ solving the equation
	\begin{align} \label{eq-MAstabtongquat0}
	(\ddc u_j+ \theta_j)^n= e^{u_j} |s|_h^{-2}f_j  \omega^n.
	\end{align}	
We have the following standard observation. 

\begin{lemma} \label{le-upperboundforuj}
There exists a constant $C>0$ such that $u_j \le C$ for every $j$. 
\end{lemma}

\proof
Let $c_j:= \int_X f_j\omega^n$ which converges to $c:= \int_X f \omega^n>0$ as $j \to \infty$. Let $A>0$ be such that $\theta_j\leq A\omega$ for all $j$.
	Since the set $\{u\in\PSH(X,A\omega):\sup_Xu=0 \}$ is compact in the $L^1$ topology,
we have \[ \int_X(u_j-\sup_Xu_j)f\omega^n \ge   - C_1\]
 for $C_1>0$ only depending on $A$ and the measure $f\omega^n$, cf.~\cite[Proposition 8.5]{GZbook}. By the monotonicity of  the sequence $(f_j)$ we have
	\begin{align*}
	    \int_Xu_j f_j\omega^n&= c_j \sup_X u_j+\int_X(u_j-\sup_Xu_j)f_j\omega^n\\
     &\geq c_j \sup_X u_j+\int_X(u_j-\sup_Xu_j)f\omega^n \ge  c_j \sup_X u_j - C_1
	\end{align*}
 Using Jensen's inequality and the fact that $|s|_h \le 1$, we have that
	\begin{align*}
	    \int_X u_j f_j\omega^n&\leq c_j\log\left( \frac{1}{c_j}\int_X e^{u_j}f_j\omega^n\right)\le c_j\log\left( \frac{1}{c_j}\int_X (\ddc u_{j}+\theta_j)^n\right)\\
     &=c_j\log\left(\frac{1}{c_j} \vol(\{\theta_j\})\right)\leq C
	\end{align*} for $C>0$ independent of $j$. This provides a uniform upper bound for $u_j$.
\endproof
	
	Our main result in this section is Proposition \ref{pro-loweerbounduepsiKE} providing a uniform lower bound for $u_j$. We use an idea in \cite{DDG-family}.

	Let $\gamma \in (0,1/2]$ be a fixed constant. 
	For a closed smooth form $\eta$, we denote by $\{\eta\}$ the cohomology class of $\eta$, and if $E$ is a divisor, we denote by $\{E\}$ the cohomology class of $E$.  
	
	Let $\mathcal{I}_D$ be the ideal sheaf defining $D$.   
   Since $\{\theta\}$ is big, we can find a K\"ahler current $T= \ddc \rho+ \theta$ in $\alpha$ such that $T\geq \delta\omega$ for some $\delta>0$. It follows from \cite{Demailly_regula_11current} that we can further assume that
$T$ has {\em analytic singularities}, i.e., there exists $c>0$ such that locally on $X$
\[\rho=\frac{c}{2}\log\sum_{\ell=1}^N|f_\ell|^2+ v \]
 where $v$ is smooth and $f_1,\ldots,f_N$ are local holomorphic functions.  The coherent ideal
sheaf $\mathcal{I}_T$ locally generated by these functions is then globally defined.

	\begin{lemma} \label{le-fujita}   There exists $\pi= \pi_\gamma: \widehat X \to X$ a finite composition of blowups with smooth centers, a smooth K\"ahler form $\widehat \omega$ and an effective divisor $E$ on $\widehat X$ so that the following conditions are satisfied:
		
		(i) the total transform $\pi^* (\mathcal{I}_D \cdot \mathcal{I}_T)$ of $\mathcal{I}_D \cdot \mathcal{I}_T$ is generated by  a divisor $\tilde{D}$, and $\supp \tilde{D} \cup \supp E$ is of simple normal crossings,
		
		(ii)
		$$\pi^* \{\theta\}= \{\widehat \omega\}+ \{E\},$$
		and
		\begin{align}\label{ine-Fujita}
		\int_{\widehat X} \widehat \omega^n \ge \vol(\{\theta\})- \gamma.
		\end{align}
	\end{lemma}
	
	\proof By \cite[Proposition 1.19]{BEGZ}, there exists $\pi': \widehat X' \to X$ a finite composition of blowups with smooth centers so that 
	$$\pi^* \{\theta\}= \{\widehat \omega'\}+ \{E'\},$$
	where $\widehat \omega'$ is a K\"ahler form on $\widehat X'$ and $E'$ is an effective divisor, and 
	\begin{align}\label{ine-Fujita2}
	\int_{\widehat X} \widehat \omega'^n \ge \vol(\{\theta\})- \gamma/2.
	\end{align}
	Let $E'$ be the exceptional divisor of $\pi'$ and $\mathcal{I}_{E'}$ be the ideal sheaf generated by $E'$. 
	By performing some more blowups to principalize $\mathcal{I}_{E} \cdot \pi'^* \mathcal{I}_D \cdot \mathcal{I}_T$, we obtain a composition of blowups  $\pi'': \widehat X \to \widehat X'$ so that if $\pi:= \pi'' \circ \pi$, then $\pi^* \mathcal{I}_D \cdot \mathcal{I}_T$ is given by some divisor  $\tilde{D}$ with simple normal crossings support and the support of this divisor also has a simple normal crossing with exceptional divisors. Repeating the arguments from \cite[Proposition 1.19]{BEGZ} shows that $\pi$ satisfies the required properties.  
	\endproof
	We denote by $\widehat D$ the support of the divisor generating $\pi^* \mathcal{I}_D$. Hence, $\widehat D$ is of simple normal crossings.  
	Let $E$ be the exceptional divisor of $\pi$ defined in Lemma~\ref{le-fujita}. 	Let 
	$$ \mu_j:= f_j |s|_h^{-2} \omega^n.$$
	Since  $u_j \in \mathcal{E}(X, \theta_j)$ is the solution of the equation
	\begin{align*}
	(\ddc u_j+ \theta_j)^n = e^{u_j} \mu_j,
	\end{align*}
	we get 
	$$(\ddc (u_j \circ \pi)+ \pi^* \theta_j )^n= e^{u_j \circ \pi} \pi^* \mu_j.$$
	Let $\widehat s$ be a section in $\mathcal{O}_{\widehat X}(\widehat D)$ defining $\widehat D$. Let $\widehat h$ be a smooth Hermitian metric on $\mathcal{O}_{\widehat X}(\widehat D)$ whose curvature is a (1,1) form $\Theta_{\widehat h}(\widehat D)$. We normalize $\widehat h$ so that $|\widehat s|_{\widehat h}$ is very small (to be made precise later). Denote $\widehat \mu_j: = \pi^* \mu_j$, and  $\widehat \mu: = \pi^* \mu$ and 
	$$\widehat f_j:= \pi^* \mu_j/ (|\widehat s|_{\widehat h}^{-2} \widehat \omega^n), \quad  \widehat f:= \pi^* \mu/ (|\widehat s|_{\widehat h}^{-2} \widehat \omega^n).$$
 With these notations, we have $\widehat\mu_j=\widehat f_j|\widehat s|^{-2}_{\widehat h}\widehat\omega^n$, $\widehat\mu=\widehat f|\widehat s|^{-2}_{\widehat h}\widehat\omega^n$.
 
	\begin{lemma}\label{le-tinhpullbackvolume}  We have that $\widehat f_j,  \widehat f$ are continuous functions and   $\widehat f_j$ increases pointwise to $\widehat f$.
	\end{lemma}
	
	\proof It suffices to work locally. Let $(\widehat z_1, \ldots, \widehat z_n)$ be a local coordinate system around a point $\widehat a$ in $\widehat X$ and $(z_1, \ldots, z_n)$ be a local coordinate system around $\pi(\widehat a)$ such that 
    \begin{enumerate}[label=\roman*., itemsep=0pt, topsep=0pt]
\item[(i)] $s(z)= z_1 \cdots z_k$ near $\pi(\widehat a)$,
\item[(ii)]  $\widehat D$ is given near $\widehat a$ by the equation $\widehat z_1 \cdots \widehat z_m =0$,
\item[(iii)] $z_j \circ \pi= \widehat z_1^{r_{1j}} \cdots \widehat z_m^{r_{mj}}$ for $1 \le j \le k$, and  some positive integers $r_{1j}, \ldots, r_{mj}$.
\end{enumerate}
Note that in (iii), since $\pi$ is a finite composition of blow-ups with smooth centers, the integers are independent of the choice of local charts; see, e.g.,~\cite[page 184]{Griffiths-Harris}.

	It follows that 
	$$\pi^*(dz_j \wedge d \bar z_j)= \sum_{1 \le p,q \le m} a_{pq} \widehat z_p^{-1} \overline{\widehat{z_q}}^{-1} |z_j \circ \pi|^2 d \widehat z_p \wedge d \overline{ \widehat{z_q}},$$
	for some constants $a_{pq}$ and $1  \le j \le k$. Hence 
	$$\pi^* \omega^n =  a \prod_{j=1}^k |z_j \circ \pi|^2  \prod_{j=1}^m |\widehat z_j|^{-2} \widehat \omega^n= a'  |s\circ \pi|_h^2 |\widehat s|_{\widehat h}^{-2} \widehat\omega^n,$$
	where $a$ and $a'$ are  smooth functions. 
 We deduce that $\widehat f=a'(f\circ \pi)$, $\widehat f_j=a'(f_j\circ\pi)$.
Therefore, $\widehat f_j,\widehat f$ are continuous functions, and $\widehat f_j$ increases pointwise to $\widehat f$ as desired.   
	\endproof
	Since $\pi^*\theta-\widehat{\omega}$ is cohomologous to $[E]$, we get 
	\begin{equation}\label{eq: phiE}
	    \ddc\varphi_E+\pi^*\theta-\widehat{\omega}=[E]
	\end{equation}
	for some negative $(\pi^*\theta-\widehat{\omega})$-psh function $\varphi_E$.
	Set 
	$$\widehat{\phi}:=-(n+2)\log(-\log|\widehat s|^2_{\widehat h} )+\varphi_E.$$
	We compute
	\begin{align}\label{eq: ddcpsihat}
	\ddc\widehat{\phi}&=-\frac{(n+2)\Theta_{\widehat h}(\widehat D)}{(-\log|\widehat s|^2_{\widehat h})}+(n+2)\frac{d  \log |\widehat{s}|^2_{\widehat h}\wedge d^c \log |\widehat{s}|^2_{\widehat h}}{(-\log|\widehat s|^2_{\widehat h})^2}+ [E]+\widehat{\omega}-\pi^*\theta.
    \end{align} 
    For simplicity, we will omit the subscript $\widehat h$. 
	Consequently, we have
	\begin{align}\label{eq: ddcpsihatthem}
	\ddc\widehat{\phi}+ \pi^* \theta  \ge 3\widehat{\omega}/4
	\end{align}
	because $|\widehat s|$ is very small. Moreover, there holds
	\[ e^{\widehat{\phi}}\widehat{\mu}\leq \frac{\widehat f\widehat{\omega}^n}{|\widehat{s}|^2(-\log|\widehat s|^2)^{n+2}}\lesssim g \widehat{\omega}^n, \quad\text{where}\; g= \frac{1}{|\widehat{s}|^2(-\log|\widehat s|^2)^{n+2}}.\] 
    By the argument in~\cite[Lemma 4.6]{DDG-family} the density $g$ satisfies 
	\[\int _{\widehat X} g|\log g|^{n+1/2}\widehat{\omega}^n<\infty.\] Indeed, $\log(-\log|\widehat s|^2)$ can be absorbed by the more singular one $-\log |\widehat s|^2$. The integral is dominated by \[C\int_{\widehat X}\frac{1}{|\widehat s|^2(-\log|\widehat s|)^{n+2}} (-\log |\widehat s|)^{n+1/2}{\widehat\omega^n}\lesssim\int_{\widehat X}\frac{\widehat\omega^n}{|\widehat s|^2(-\log|\widehat s|)^{3/2}}. \]
 Since $\widehat D$ has snc support, and $|\widehat s|<1$, it suffices to check that on the unit polydisc $\Delta^n\subset\mathbb{C}^n$,
 the integral \[ \int_{\Delta^n}\frac{dV_{\mathbb{C}^n}}{\prod_{i=1}^k|z_i|^2(-\log|z_i|)^{3/2}}\] converges.
 Using polar coordinates,
the conclusion follows from $\int_0^{1}\frac{dr}{r(-\log r)^{3/2}}<\infty$.
 
 By \cite{Kolodziej05} or ~\cite[Theorem 1.5]{DDG-family}, we can find $\varphi$ a bounded $\widehat \omega/2$-psh function satisfying
	\begin{align}\label{eq: phi}
	(\ddc \varphi+ \widehat \omega/2)^n= e^{\varphi+\widehat \phi} \widehat \mu.
	\end{align}
	We note that $\varphi$ is bounded globally. Define
	\begin{equation}\label{eq: psihat}
	    \widehat\psi_\gamma:= \varphi+ \widehat \phi 
	\end{equation}
	which is a $\pi^* \theta_j$-psh function for $j \ge j_\gamma$ large enough since it follows from~\eqref{eq: ddcpsihatthem} that
	\begin{align*}
	\ddc\widehat\psi_\gamma+\pi^* \theta&\ge  \ddc \varphi +3\widehat\omega/4\geq   \widehat\omega/4,
	\end{align*}
	and $\pi^* \theta_j- \pi^* \theta \le \widehat \omega/4$ for $j$ big enough (depending on $\gamma$). 
	
	\begin{proposition}\label{pro-loweerbounduepsiKE} For any $\gamma>0$,
 we have
		\begin{equation}\label{ine-GDG}
		u_j\circ \pi \ge \widehat \psi_\gamma,
		\end{equation}
		for $j\geq j_\gamma$. Moreover, there exists a smooth Hermitian metric $\widehat h$ (depending on $\gamma$) on $\mathcal{O}_{\widehat X}(\widehat D)$ such that
		$$\int_{\widehat X} (\ddc \widehat \psi_\gamma+\pi^* \theta_j)^n \ge  \vol(\{\theta_j\}) - 2\gamma$$
		for $j$ big enough (depending on $\gamma$).
	\end{proposition}
	
	\proof
	Observe that 
	$$(\ddc \widehat \psi+ \pi^* \theta_j)^n \ge (\ddc \varphi+ \widehat \omega/2)^n = e^{\varphi+\widehat \phi} \widehat \mu.$$
	The domination principle (note that $\widehat \psi$ might not be of full mass, but $u_j$ is so) yields the first desired assertion. Indeed, from the above inequality, we have
	\[e^{-\widehat\psi} (\ddc \widehat \psi+ \pi^* \theta_j)^n\geq e^{-u_j\circ \pi}(\ddc  u_j\circ \pi+ \pi^* \theta_j)^n.\]
Therefore, $u_j\circ \pi \ge \widehat \psi$ on $\widehat X$ by  Proposition~\ref{prop: standarddomination}.
	
We now check the second desired assertion. Put $v_1:= -(n+2)\log(-\log |\widehat s|_{\widehat h}^2)$, and $v_2:=  \log |\widehat s|_{\widehat h}^2$. We choose $\widehat h$ so that $|\widehat s|_{\widehat h}^2$ is so small that $v_2$ is $\gamma\widehat \omega$-psh function with $\gamma\in(0,1/2]$, and $v_1 \ge -(-v_2)^{1/2}-C$ by the elementary inequality $(n+2)\log y+O(1)\leq y^{1/2}$ for $y>0$. 
It follows that $ -(-v_2)^{1/2} \in \mathcal{E}(\widehat X, \gamma \widehat \omega)$ (by \cite[Example 2.7]{ComanGuedjZeriahi-DMA}), hence, $v_1 \in \mathcal{E}(\widehat X,\gamma \widehat \omega)$, cf.~\cite[Proposition 1.6]{GZ-weighted}.  Direct computations show
\begin{align*}
  \ddc \widehat\psi_\gamma + \pi^* \theta&= \ddc v_1+ \ddc \varphi_E+ \ddc \varphi+ \pi^* \theta\\
  &\ge (\ddc v_1+ \gamma \widehat \omega) + (\ddc \varphi+ \widehat \omega/2)+ (1- \gamma -1/2) \widehat \omega
\end{align*}
using the equality~\eqref{eq: phiE} in the last inequality.
Let $R$ be the right-hand side of the last inequality. Since $v_1 \in \mathcal{E}(\widehat X, \gamma \widehat \omega)$ and $\varphi \in \mathcal{E}(\widehat X, \widehat\omega/2)$ (defined in~\eqref{eq: phi}), one infers that 
$$\int_{\widehat X} R^n = \int_{\widehat X} \widehat \omega^n \ge \vol(\{\theta\})-\gamma.$$
Consequently
$$\int_{\widehat X} (\ddc \widehat\psi_\gamma + \pi^* \theta)^n \ge \vol(\{\theta\})- \gamma.$$
Since $\theta_j$ is close to $\theta$ in the sup norm, we deduce that 
$$\int_{\widehat X} (\ddc \widehat\psi_\gamma + \pi^* \theta_j)^n \ge \vol(\{\theta_j\})- 2\gamma$$
for $j\geq j_\gamma$ big enough. We note that $\pi^*\theta_j\to\pi^*\theta$ in the $\mathcal{C}^0$ topology as $j\to\infty$.
This finishes the proof.
	\endproof 
	
	
	\begin{corollary} \label{cor-subsolution} Let $u$ be a $L^1$ limit of $(u_j)_j$ as $j \to \infty$. Then $u \in \mathcal{E}(X, \theta)$ and 
\begin{align} \label{ine-reviseduliminf}		
		(\ddc u+ \theta)^n \ge e^u f |s|_h^{-2} \omega^n.
		\end{align}
	\end{corollary} 
	
	\proof Without loss of generality, we can assume that $u_j \to u$ in $L^1$. Let $\epsilon>0$ be a constant, and $\theta'_\epsilon:= \theta+ \epsilon \omega$. Let $j_\epsilon \in \N$ be such that $\theta_j \le \theta'_\epsilon$ for every $j \ge j_\epsilon$. 
Fixing $\gamma\in (0,1/2]$, we consider the modification $\pi:\widehat X\to X$ as in Lemma~\ref{le-fujita}, and $\widehat\psi_\gamma$ is defined in~\eqref{eq: psihat}.
 Let $\psi_\gamma:= {\pi}_* \widehat \psi_\gamma$ which is a $\theta_j$-psh function for $j$ large enough.
 It follows that $\psi_\gamma$ is $\theta'_\epsilon$-psh and $\ddc \psi_\gamma + \theta'_\epsilon \ge \ddc \psi_\gamma + \theta_j$ for $j \ge j_\epsilon$.  
 Since $\pi$ is a bimeromorphic map, one obtains
$$\int_X (\ddc \psi_\gamma+ \theta_j)^n= \int_{\widehat X} (\ddc \widehat \psi_\gamma+ \pi^* \theta_j)^n.$$ 	
	 By this and Proposition~\ref{pro-loweerbounduepsiKE}, we get $u \ge \psi_\gamma$ and 
	$$\int_X (\ddc \psi_\gamma+ \theta_j)^n \ge \vol(\{\theta\}) -2\gamma.$$ 
	This, combined with the monotonicity of non-pluripolar products (see \cite{Lu-Darvas-DiNezza-mono}), gives
\begin{align*}	
	\int_X (\ddc u+ \theta'_\epsilon)^n &\ge \int_X (\ddc \psi_\gamma+\theta'_\epsilon)^n \\
	& \ge \int_X (\ddc \psi_\gamma+ \theta_j)^n  \ge \vol(\{\theta\})- 2 \gamma
\end{align*}	
	for every constant $\gamma >0$. Letting $\gamma \to 0$ gives
\begin{align} \label{inetraloirefereupsigamma}
\int_X (\ddc u+ \theta'_\epsilon)^n \ge  \vol(\{\theta\}).
\end{align}	
 We underline here that since $u$ is independent of $\gamma$, we can take $\gamma \to 0$ in the above inequality. Taking $\epsilon \to 0$, we obtain $u \in \mathcal{E}(X, \theta)$.    
	 
It remains to check (\ref{ine-reviseduliminf}), whose proof is quite standard; see, e.g.,~\cite[Lemma 2.8]{Lu-Darvas-DiNezza-logconcave}. Since some slight modifications are needed (due to the fact that $\theta_j$ is not in the same cohomology class as $\theta$), we present the details for the reader's convenience. Let $u'_j:= (\sup_{k \ge j} u_k)^*$ which are $\theta'_\epsilon$-psh for $j \ge j_\epsilon$. Note that $u'_j$ decreases to $u$. By extracting a subsequence, if necessary, we can assume that 
	$$(\ddc u'_j+ \theta'_\epsilon)^n \to \nu$$
	weakly as $j \to \infty$. Since $u'_j$ decreases to $u$, using \cite[Theorem 2.3]{Lu-Darvas-DiNezza-mono}, we have 
$$\nu \ge (\ddc u+ \theta)^n.$$	 	
Using this and  the fact that $u \in \mathcal{E}(X, \theta)$, we obtain that 
	\begin{align}\label{ine-sosanhsubsolution} \nonumber
	\|\nu- (\ddc u+ \theta)^n\| &\le \nu(X)- \int_X (\ddc u+ \theta)^n \\
	&\le \vol(\{\theta'_\epsilon\})- \vol(\{\theta\}) \lesssim \epsilon,
	\end{align}
	where $\|\nu'\|$ denotes the mass norm of a (signed) measure of $\nu'$.
	On the other hand, for every continuous function $g$ with compact support in $X \backslash W$ ($W$ is the union of $D$ and the non-K\"ahler locus of $\{\theta\}$), we have 
	$$g (\ddc u'_j + \theta'_\epsilon)^n \ge g e^{\inf_{k \ge j} u_j} f_j |s|^{-2}_h \omega^n$$
	which converges to $g e^{u} f |s|^{-2}_h \omega^n$ as $j \to \infty$ by the Lebesgue dominated convergence theorem (note that $\inf_{k \ge j} u_j$ converges pointwise almost everywhere to $u$). Hence letting $j \to \infty$ gives
	$$g \nu \ge g e^u f |s|^{-2}_h \omega^n.$$
	Combining this with (\ref{ine-sosanhsubsolution}) implies
	$$g (\ddc u+ \theta)^n \ge  g e^u f |s|^{-2}_h \omega^n - C\|g\|_{L^\infty} \epsilon$$
	for every constant $\epsilon>0$. Letting $\epsilon \to 0$, we obtain the desired inequality. This finishes the proof.
	\endproof
	
	Let $g$ be a function on $X$ and let $\eta$ be a smooth closed $(1,1)$-form in a pseudoeffective class. If there is a function $v \in \PSH(X, \eta)$ with $v \le g$, then we define
	$$P_\eta(g):= \big(\sup\{v \in \PSH(X, \eta): v \le g\}\big)^*$$
which is a well-defined $\eta$-psh function. 
	Let $\epsilon>0$ be a small constant. Let $\theta'_\epsilon:= \theta+ \epsilon \omega$. Let $j_\epsilon \in \N$ be such that $\theta_j \le \theta'_\epsilon$ for every $j \ge j_\epsilon$.  Hence $u_j$ is $\theta_\epsilon'$-psh for $j \ge  j_{\epsilon}$. 
	
	\begin{corollary} \label{cor-envelopchanduoi}  Let 
		$$u''_{j,\epsilon}:= \lim_{k\to \infty} P_{\theta'_\epsilon}( \min \{u_j, \ldots, u_k\}),$$
		for $j \ge j_\epsilon$.  Then $u''_{j,\epsilon}$ is  a well-defined $\theta'_\epsilon$-psh function and increases to some  $\theta'_\epsilon$-psh function $u''_\epsilon$ as $j\to \infty$ satisfying:
		
		(i) $u''_\epsilon \le u$, where $u$ is a $L^1$-limit of the sequence $(u_j)_j$ as $j\to \infty$,
		
		(ii) $u''_\epsilon$ decreases to some $\theta$-psh function $u'' \in \mathcal{E}(X, \theta)$ as $\epsilon$ decreases to $0$. 
	\end{corollary}
	
	\proof Since $u_j \ge \psi_\gamma$ for $j \ge j_\gamma$ big enough (depending on $\gamma$), we see that $u''_{j,\epsilon}$ is well-defined for $j \ge j_\gamma$ because $\psi_\gamma$ is in the envelope defining $P_{\theta'_\epsilon}(\min \{u_j, \ldots, u_k\})$ for $k \ge j$ (and $(u_j)_j$ is bounded from above uniformly in $j$, see Lemma \ref{le-upperboundforuj}). The fact that $u''_{j,\epsilon}$ is increasing (for $\epsilon$ fixed) is clear. By the definition of the envelope $P_{\theta'_\epsilon}$, one sees that $u''_\epsilon$ is decreasing as $\epsilon\searrow 0$. Since $u''_\epsilon \ge \psi_\gamma$, we see that $u'':= \lim_{\epsilon \to 0} u''_\epsilon \ge \psi_\gamma$
	for every $\gamma \in (0, 1/2]$. Using this and arguing as at the beginning of the proof of Corollary \ref{cor-subsolution} (see (\ref{inetraloirefereupsigamma})), we obtain  
	\begin{align*} 
\int_X (\ddc u+ \theta'_\epsilon)^n \ge  \vol(\{\theta\})-2\gamma
\end{align*}	
for 	every $\gamma \in (0, 1/2]$ and every constant $\epsilon>0$, where $\theta'_\epsilon:= \theta+ \epsilon \omega$. Letting $\gamma \to 0$ and then $\epsilon \to 0$ gives
	$u'' \in \mathcal{E}(X, \theta)$. 
	\endproof

	\section{Almost complete K\"ahler--Einstein metrics}

	In this section, we prove Theorem~\ref{the-TianYaumetric}.
	We begin by recalling the following definition for almost complete metrics. In what follows, we interchangeably use the terms K\"ahler forms and K\"ahler metrics. We start with a slightly more general version of Yau's Schwarz lemma for volume forms (\cite{Yau-schwarz} and also \cite{Tosatti-schwarz}). 
	
	\begin{theorem} \label{the-yauschwarz-sing} Let $(M,\eta_1)$ be a complete K\"ahler manifold with Ricci curvature bounded from below and the scalar curvature bounded from below by a constant $-\frac{nK_1}{2\pi}$. Let $N$ be a complex manifold with $\dim N= \dim M$ and $\eta_2$ be a closed positive $(1,1)$-current on $N$ having a well-defined Ricci curvature. Let $f: M \to N$ be a non-degenerate holomorphic map. Assume that the following conditions are satisfied:
		
		(i) $2 \pi \ric \eta_2 \le - K_2 \eta_2$ as currents for some constant $K_2>0$,
		
		(ii) There exists a closed subset $E$ in $N$ such that $f^{-1}(E)$ is of zero Lebesgue measure in $M$, and  $\eta_2$ is a smooth K\"ahler form on $N \backslash E$, and $\eta_2^n$ extends to a smooth form on $N$. 
		
		Then we have $K_1 \ge 0$, and 
		$$f^* \eta_2^n \le \bigg(\frac{K_1}{K_2}\bigg)^n \eta_1^n$$ 
	\end{theorem}

We note that in (i) there is the factor $2\pi$ because our $\ric \eta_2$ (see Introduction) is, in fact, equal to the $\frac{1}{2\pi}$ times the usual Ricci curvature of the metric $\eta_2$. 
	
	\proof  
	The proof is almost identical to that of~\cite[Theorem. 1.2]{Tosatti-schwarz} (which goes back to~\cite{Yau-schwarz}).
	Set
	 $$u:=\frac{f^*\eta_2^n}{\eta_1^n} \cdot$$
	 Since $\eta_2^n$ is a smooth form on $N$, we see that $u$ is a well-defined nonnegative smooth function on $M$. 
	Let $\Delta$ denote the Laplacian with respect to $\eta_1$. 
	At any point $x\in M\backslash f^{-1}(E)$, we compute (recall $\ddc:= \frac{i}{\pi} \partial \bar \partial$)
	\begin{align*}
	\Delta u&=u\Delta\log u+\frac{|\nabla u|^2_{\eta_1}}{u}\\
	&= \frac{2n  \pi  \ddc \log u \wedge \eta_1^{n-1}}{\eta_1^n}+\frac{|\nabla u|^2_{\eta_1}}{u}\\  
	&=2 \big(-{\rm tr}_{\eta_1}f^*(2\pi\ric\eta_2)+{\rm tr}_{\eta_1}(2 \pi\ric\eta_1)\big)u+\frac{|\nabla u|^2_{\eta_1}}{u}\\
	&\geq 2 K_2{\rm tr}_{\eta_1}f^*\eta_2\cdot u-2 nK_1\cdot u\\
	&\geq 2nK_2 u^{1+\frac{1}{n}}-2nK_1 u,
	\end{align*} using the arithmetic-geometric mean inequality. Since both sides of the last inequality are continuous on $M$, we infer
	$$\Delta u \ge 2n K_2 u^{1+ \frac{1}{n}}- 2n K_1 u$$	
	on $M$. Now, it suffices to repeat arguments in \cite[Pages 200-201]{Yau-schwarz} to obtain the desired estimate.
	\endproof

	Recall that a closed positive $(1,1)$-current on a complex manifold $M$ is said to be \emph{a singular K\"ahler--Einstein metric of negative Ricci curvature} on $M$ if $\eta$ has a well-defined Ricci curvature and $\ric \eta= - \eta$.

	\begin{definition} \label{def-almost-complete} (\cite{Tian-Yau-complete-ssurvey} and also \cite{Cheng-Yau})
		Let $M$ be a K\"ahler manifold and let $\eta$ be a singular K\"ahler--Einstein metric of negative Ricci curvature on $M$. 
		We say that $\eta$ is \emph{almost complete} if there exists a proper analytic subset $E \subset M$, and  a sequence of complete smooth K\"ahler metrics $\omega_j$ on $M$ and a sequence of positive numbers $t_j$ converging to 1 such that the following conditions are fulfilled:
		
		(i) $\ric \omega_j \ge - t_j \omega_j$ for every $j \in \N$, 
		
		(ii) $\eta$ is the limit  of $\omega_j$ in the local $\mathcal{C}^\infty$ topology on $M \backslash E$ as $j \to \infty,$ 
		
		(iii) $\eta$ is a smooth K\"ahler form on $M \backslash E$. 
	\end{definition}
	
	We stress that one should consider an almost-complete metric on $M$ rather than a smooth metric on $M \backslash E$. The key is the following property of almost-complete metrics. 
	
	\begin{lemma}[{\cite{Tian-Yau-complete-ssurvey}}]\label{le-tinhduynhatganday} Let $M$ be a K\"ahler manifold. Then there exists at most one almost complete K\"ahler--Einstein metric $\eta$ on $M$ with $\ric \eta= - \eta$. 
	\end{lemma}
	
	\proof We reproduce the proof from \cite{Tian-Yau-complete-ssurvey} for readers' convenience. Suppose that $\eta_1, \eta_2$ are two almost complete K\"ahler metrics on $M$ with $\ric \eta_k= - \eta_k$ for $k=1,2$. By definition, there are sequences of complete K\"ahler metrics $\omega_{jk}$ and positive real numbers $t_{jk} $ converging to $1$ such that 
	$$\ric \omega_{jk}  \ge - t_{jk} \omega_{jk}$$
	and $\omega_{jk}$ converges to $\eta_k$ smoothly pointwise.  By Theorem \ref{the-yauschwarz-sing} applied to the identity map of $M$, we get 
	$$\eta_1^n \le t_{j2}^n \omega_{j2}^n, \quad \eta_2^n \le t_{j1}^n \omega_{j2}^n$$
	for every $j \in \N$. Letting $j \to \infty$ and using the pointwise convergence of $\omega_{jk}$ to $\eta_k$ give $\eta^n_1= \eta^n_2$. Thus 
	$$\eta_1= -\ric \eta_1= \ddc \log \eta_1^n=  \ddc \log \eta_2^n = - \ric \eta_2= \eta_2.$$
	This finishes the proof.   
	\endproof
	
	\begin{lemma}\label{le-uniform-Cinfintynorm}
		Let $D$ be a simple normal crossing divisor such that $K_X+ D$ is big and nef. Let $s$ be a section of $\mathcal{O}_X(D)$ defining $D$ and $h$ be a smooth Hermitian metric on $\mathcal{O}_X(D)$.
		Let $\epsilon \in [0,1)$ and $f$ be a smooth function. Let  $\theta$  be a smooth representative of big cohomology class $c_1(K_X+D)$. Let $v_\epsilon \in \mathcal{E}(X, \theta+\epsilon\omega)$ be the unique solution of the equation
		\begin{equation}\label{eq: cmae2}
		(\ddc v_\epsilon+\theta+\epsilon\omega)^n=e^{v_\epsilon+f}|s|_h^{-2}\omega^n.
		\end{equation}  
		Let $E$ be the non-K\"ahler locus of $K_X+D$.  Then the $\mathcal{C}^k$-norm of $v_\epsilon$ on compact subsets in $X \backslash (D \cup E)$ is uniformly bounded in $\epsilon$ for every $k$. 
	\end{lemma}
	
	\proof The domination principle (see Proposition~\ref{prop: standarddomination}) ensures that $v_\epsilon$ decreases to $v$ as $\epsilon\searrow 0$. 
	It follows that $v_\epsilon\leq C$ for some uniform constant $C>0$. 
	Theorem~\ref{the-stability} shows that $v_\epsilon$ converges weakly toward $v$ as $\epsilon\to 0^+$. Moreover, $v$ belongs to $\mathcal{E}(X,\theta)$ and satisfies
	\[(\ddc v+\theta)^n =e^{v+f}|s|^{-2}_h\omega^n.\]

	We set $X_0:=X\setminus (D\cup E)$.
	We are going to prove that the family $(v_\epsilon)$ is precompact in $\mathcal{C}^\infty_{\rm loc}(X_0)$, where $E$ is the non-ample locus of $K_X+D$. This amounts to establishing $\mathcal{C}^k_{\rm loc}(X_0)$ estimates for all $k\in \mathbb{N}$ thanks to the Arzela--Ascoli theorem.
	According to the Evans--Krylov theory and Schauder interior estimates (the so-called bootstrapping arguments for elliptic PDEs), obtaining a local $L^\infty$ and Laplacian estimate on $X_0$ suffices. To do this,  we repeat the estimates derived  in~\cite[Sect. 5]{DiNezza-Lu-quasiprojective} with minor modifications.

	Since $\theta$ represents a big cohomology class, we can find a $\theta$-psh function $\rho$ with analytic singularities such that $\rho\to -\infty$ near $E$ and $\ddc\rho+\theta\geq 2\delta\omega$ for some $\delta>0$. We fix $c>0$ so small that $c\ddc\log|s|^2_h+\delta\omega\ge 0$.
	It follows from~\cite[Theorem. 3.1, Step 1]{dang22continuity} that \[v\geq c\log|s|^2_h+\rho-C \]
	for $C>0$ depending on $X$, $\omega$, $n$, $c$, $\delta$ and an upper bound for $\int_Xe^{-2P_\omega(c^{-1}(v-V_\theta))}\omega^n$. Therefore, 
	for every $\epsilon\in (0,1]$,
	\begin{equation}
	\label{eq: lowerbound-v}
	v_\epsilon\geq c\log|s|_h^2+\rho-C.
	\end{equation}
	
	We are now in a position to establish the local Laplacian estimate. 
	We recall Siu--Yau's inequality (cf.~\cite{Siu-Hermitia-Einstein-metric}): Let $\tau$ and $\tau'$ are two K\"ahler
	forms on a complex manifold and let $F$ be
	defined by $\tau'^n=e^F\tau^n$. If the holomorphic bisectional curvature of $\tau$ is bounded from below by some constant $B>0$, then
	\begin{align*}
	\Delta_{\tau'}\log{\rm tr}_\tau(\tau')\geq \frac{\Delta_\tau F}{{\rm tr}_\tau(\tau')}-B{\rm tr}_{\tau'}(\tau).
	\end{align*}

	We are going to apply $\tau=\omega$ and $\tau'=\omega_\epsilon':=\ddc v_\epsilon+\theta+\epsilon\omega$. We observe that the holomorphic bisectional curvature of $\omega$ is obviously bounded on $X$ by a constant $B>0$, hence we obtain the following inequality
	\begin{equation}\label{eq: lap}
	\Delta_{\omega_\epsilon'}\log {\rm tr}_{\omega}(\omega_\epsilon' )\geq
	\frac{\Delta_{\omega}(v_\epsilon+F)}{{\rm tr}_{\omega}(\omega_\epsilon' )}-B{\rm tr}_{\omega_\epsilon'}(\omega),
	\end{equation} 
	where $F:= f- \log |s|^2_h$.
	We have $\Delta_{\omega}F\ge  - nA$ on $X \backslash D$ for some uniform constant $A>0$.
	Moreover, one sees that 
	 $$\Delta_{\omega}v_\epsilon={\rm tr}_{\omega}(\ddc v_\epsilon+ C_0\omega)-nC_0\geq -nC_0,$$ for some uniform $C_0>0$.  Combining this  with the elementary inequality ${\rm tr}_{\omega}(\omega_\epsilon'){\rm tr}_{\omega_\epsilon'}(\omega)\geq n$, we obtain
	\begin{equation}\label{eq: lap2}
	\frac{\Delta_{\omega}(v_\epsilon+F)}{{\rm tr}_{\omega}(\omega_\epsilon')}\geq -(A+C_0){\rm tr}_{\omega_\epsilon'}(\omega).
	\end{equation}
	Plugging~\eqref{eq: lap2} into~\eqref{eq: lap}, we thus obtain
	\begin{equation*}
	\Delta_{\omega_\epsilon'}\log{\rm tr}_{\omega}(\omega_\epsilon')\geq -C_1{\rm tr}_{\omega_\epsilon'}(\omega)
	\end{equation*} for $C_1=B+A+C_0$. We set $w_\epsilon:=v_\epsilon-c\log|s|^2_h-\rho$. 
	By~\eqref{eq: lowerbound-v}, we have that $w_\epsilon$ is bounded from below.
	We set \[\omega_\epsilon:=\ddc c\log|s|^2_h+\ddc\rho+\theta+\epsilon\omega.\]
	Since we have chosen $c>0$ so small that $\ddc c\log|s|^2_h+\delta\omega\geq 0$, hence $\omega_\epsilon\geq \delta\omega$ for every $\epsilon$.
	We observe that
	$\Delta_{\omega_\epsilon'}w_\epsilon=n-{\rm tr}_{\omega_\epsilon'}\omega_\epsilon\leq n-\delta{\rm tr}_{\omega_\epsilon'}\omega$.
	Therefore, we obtain \begin{equation}\label{eq: max-prin}
	\Delta_{\omega_\epsilon'}(\log{\rm tr}_{\omega}(\omega_\epsilon')-(C_1\delta^{-1}+1)w_\epsilon)\geq {\rm tr}_{\omega_\epsilon'}(\omega)-n(C_1\delta^{-1}+1).
	\end{equation}  
	We are now in a position to apply the maximum principle. Indeed, if we put $C=C_1\delta^{-1}$ then since $w_\epsilon$ tends to $+\infty$ near $D\cup E$
	the function $H:=\log{\rm tr}_{\omega}(\omega_\epsilon')-(C+1)w_\epsilon$ therefore attains its maximum at some point $x_0\in X_0$ (depending on $\epsilon$). At this point, the inequality~\eqref{eq: max-prin} combined with the maximum principle yield ${\rm tr}_{\omega_\epsilon'}(\omega_\epsilon)(x_0)\leq n(C+1)$. Using the elementary inequality ${\rm tr}_{\tau}(\tau')\leq n\left(\frac{\tau'^n}{\tau^n} \right)({\rm tr}_{\tau'}(\tau))^{n-1}$
	for any two K\"ahler forms $\tau$, $\tau'$, one gets 
	\begin{align*}
	\log{\rm tr}_{\omega}(\omega_\epsilon')&\leq \log{\rm tr}_{\omega}(\omega_\epsilon')(x_0)+(C+1)(w_\epsilon-w_\epsilon(x_0))\\
	&\leq (n-1)\log n(C+1)+\log n+ w_\epsilon(x_0)+F_\epsilon(x_0)+(C+1)(w_\epsilon-w_\epsilon(x_0))\\
	&\leq C_3+(C+1)w_\epsilon
	\end{align*} since $w_\epsilon$  is uniformly bounded from below by~\eqref{eq: lowerbound-v}. This implies that $H$ is uniformly bounded from above, hence
	\[ {\rm tr}_{\omega}(\omega_\epsilon')\leq C_3 e^{-(C+1)\rho}\frac{1}{|s|_h^{2c(C+1)}}\] using that $v_\epsilon$ is uniformly bounded from above. 
	Thus, we end up with uniformly (in $\epsilon$) positive constant $C>0$ such that\[ |\Delta_\omega v_\epsilon|\leq \frac{Ce^{-C\rho}}{|s|_h^{2C}} \cdot \]
	In particular, for any compact $K\subset \subset X_0$, one gets a uniform bound for $\|\Delta_\omega v_\epsilon\|_{L^\infty(K)}$. We can apply a complex version of the Evans--Krylov estimate, due to Trudinger~\cite{trudinger83fully} (cf.~\cite[Chap. 2]{Siu-Hermitia-Einstein-metric} in this context) to obtain a local $\mathcal{C}^\alpha$ estimate on the metric $\omega_\epsilon'$ for $0<\alpha<1$, i.e., $\|\Delta_\omega v_\epsilon \|_{\mathcal{C}^{2,\alpha}(K)}\leq C_K$ for a uniform constant $C_K>0$ (we also refer to~\cite{WangWu20regularity} for a new proof). From this, we eventually differentiate the equation using Schauder estimates to obtain uniform estimates  
	$$\sup_{\epsilon\in (0,1]}\| v_\epsilon\|_{\mathcal{C}^{j,\alpha}(K)}\leq C_{j,\alpha}(K)<+\infty,$$ for each $0<\alpha<1$, $j\in\mathbb{N}$, which
	guarantee that $v_\epsilon$ is relatively compact in $\mathcal{C}^\infty_{\rm loc}(X_0)$. The lemma is thus proved.
	\endproof

	

	\begin{proof}[End of the proof of Theorem \ref{the-TianYaumetric}]
		Let $s$ be a section of $\mathcal{O}_X(D)$ defining $D$, and $h$ a smooth metric on $\mathcal{O}_X(D)$ such that $|s|_h \le 1$. Let $\eta$ be the Chern form of $h$. Let $\{\omega\}$ denote the cohomology class of $\omega$.
		
		Let $\theta$ be a smooth closed form in $c_1(K_X+D)$ and $\theta_\epsilon:= \theta+ \epsilon \omega$ for $\epsilon \in [0,1]$ (hence $\theta_0= \theta$). 
		Thus  $(\theta_\epsilon)_\epsilon$ is a sequence of smooth closed forms such that $\theta_\epsilon \in c_1(K_X+D)+ \epsilon \{\omega\}$ and $\theta_\epsilon$ converges to $\theta$ in $\mathcal{C}^0$-topology as $\epsilon \to 0$. For $0 \le  \epsilon \le 1$,  let $u_\epsilon \in \mathcal{E}(X, \theta_\epsilon)$ be the solution of
		$$(\ddc u_\epsilon+ \theta_\epsilon)^n= e^{2 u_\epsilon+2F} |s|^{-2}_h \omega^n,$$
		where $F$ is a smooth function so that 
		$$\eta- \ric (\omega)-\theta= -\ddc F.$$
		We see that $\omega_\epsilon:= \ddc u_\epsilon + \theta+ \epsilon \omega$ satisfies  
		$$\ric \omega_\epsilon= - \omega_\epsilon+ \epsilon \omega + [D]$$
		and $\omega_\epsilon$ is of full Monge--Amp\`ere mass in $c_1(K_X+D)+ \epsilon \{\omega\}$ for $0 \le \epsilon \le 1$.    
		Observe that the setting in Section \ref{sect: lowerbound} applies to $(u_\epsilon)_\epsilon$. Hence,  by Corollary~\ref{cor-subsolution}, we see that if $u$ is the $L^1$ limit of a subsequence $(u_{\epsilon_j})_j$ as $j \to \infty$, then $u \in \mathcal{E}(X, \theta_0)$. Lemma \ref{le-uniform-Cinfintynorm} yields that $u_{\epsilon_j}$ also converges to $u$ in the local $\mathcal{C}^\infty$ topology in $X \backslash (D \cup E)$. Consequently, we get 
		$$(\ddc u+ \theta_0)^n =  e^{2 u+2F} |s|^{-2}_h \omega^n$$ 
		on $X \backslash (D \cup E)$.  It follows that $u=u_0$ by uniqueness (see \cite{BEGZ} or \cite{Dinew-uniqueness}). Consequently $u_\epsilon$ converges to $u_0$ as $\epsilon \to 0$ in $L^1$ and $u_\epsilon$ converges to $u_0$ in the local $\mathcal{C}^\infty$-topology in $X \backslash (D \cup E)$. Since $\omega_\epsilon$ is a complete smooth K\"ahler metric (by \cite{Kobayashi-KE}) on $X \backslash D$ and $\ric \omega_\epsilon \ge - \omega_\epsilon$, we obtain that $\omega_0$ is almost-complete, and the uniqueness of $\omega_0$ follows from Lemma  \ref{le-tinhduynhatganday}. The sought metric $\omega_D$ is exactly $\omega_0$ in the above discussion. This finishes the proof.
	\end{proof}

	We end this section by presenting an example in which the non-K\"ahler locus of $K_X+D$ is not contained in $D$. 

\begin{example} \label{ex-blowupP2} 
	Fix now a point $x_0 \in \P^2$. Let $\pi= \pi_{x_0}:X \to \P^2$ be the blowup $\P^2$ at a given point $x_0$ and let $E$ be the exceptional divisor. Let $D'= H_1+ H_2+H_3+H_4$ be a smooth divisor in $\P^2$ where $H_1, H_2,H_3,H_4$ are hyperplanes in general position, and $x_0 \in H_1 \backslash (\cup_{j=2}^4 H_j)$. Let $\pi^* D'$ be the total transform of $D'$ and let $D$ be the strict transform of $D'$. They are both simple normal crossing divisors in $X$.
 By our choice of $x_0$, we get 
	$$\pi^*D'= D+ E.$$
 Observe that $\pi^* K_{\P^2}=K_X- E$ (see, e.g.,~\cite[page 187]{Griffiths-Harris}) and $K_{\P^2}+ D'= \mathcal{O}_{\P^2}(1)$ which is a positive line bundle (see, e.g.,~\cite[page 146]{Griffiths-Harris}). We infer that 
	$$K_X+ D= \pi^*K_{\P^2}+E +\pi^*D' -E= \pi^*(K_{\P^2}+D')= \pi^* \mathcal{O}_{\P^2}(1)$$
	which is a semi-ample and big line bundle on $X$. One can see that $E$ is equal to the non-K\"ahler locus of $K_X+D$, but is not contained in $D$. 
\end{example}

	

	\section{Stability of complex Monge--Amp\`ere equations}

	The goal of this section is to prove the following result.

	\begin{theorem} \label{the-stability} Let $D$ be a simple normal crossing divisor in $X$. Let $s$ be a section of $\mathcal{O}_X(D)$ defining $D$ and $h$ be a smooth Hermitian metric on $\mathcal{O}_X(D)$ such that $|s|_h<1$. Let $\alpha$ be a big cohomology class and $\theta$ be a smooth representative in $\alpha$. Let $(f_j)_j$ be an increasing sequence of continuous nonnegative functions such that $f_1 \not \equiv 0$, and $f_j$  converges pointwise to $f \in \mathcal{C}^0(X)$ as $j \to \infty$. Let $(\theta_j)_{j\in \mathbb{N}}$ be a sequence of smooth closed $(1,1)$-forms in big cohomology classes converging to $\theta$ in the $\mathcal{C}^0$ topology as $j \to \infty$. Let $u_j \in \mathcal{E}(X, \theta_j)$ be a solution of the equation
		\begin{align} \label{eq-MAstabtongquat}
		(\ddc u_j+ \theta_j)^n= e^{u_j} |s|_h^{-2}f_j  \omega^n
		\end{align}
		for every $j \ge 0$. Then $(u_j)_{j\in \mathbb{N}}$ is a decreasing sequence in capacity and $u_j$ converges in $L^1$ to a $\theta$-psh function $u \in \mathcal{E}(X, \theta)$ satisfying
		\begin{align} \label{eq-MAstabtongquat2}
		(\ddc u+ \theta)^n= e^{u} |s|_h^{-2}f  \omega^n.
		\end{align}
		In particular if (\ref{eq-MAstabtongquat}) admits a solution in $\mathcal{E}(X, \theta_j)$, then (\ref{eq-MAstabtongquat2}) also possesses a (unique) solution in $\mathcal{E}(X, \theta)$.
	\end{theorem}

	To see the relevance of this result to  Theorem \ref{the-KE-Gue}, 
	we recall the Monge--Amp\`ere equation satisfied by $\omega_\epsilon$. 
	Let $\theta$ be a smooth closed form representing $c_1(K_X+D)$. Thus $\omega_D= \ddc u+ \theta $ where $u$ is an unknown $\theta$-psh function. Let $s$ be a section in $\mathcal{O}_X(D)$ defining $D$. Let $h$ be a smooth Hermitian metric on the latter line bundle. Let $\eta$ be the Chern form of $h$. We rescale $h$ so that $|s|_h^2 \le 1$.  Observe that 
	$$[D]= \ddc \log |s|_h + \eta,$$
	and  by $\ddc$--Lemma one gets
	$$\eta - \ric (\omega)-\theta= -\ddc F,$$
	for some smooth function $F$ on $X$. On the other hand,
	$$2\ric (\omega_D)= - \ddc \log (\omega_D^n/\omega^n)+ 2\ric (\omega).$$
	Hence since  $\ric \omega_D= - \omega_D+ [D]$, one obtains
	$$\frac{1}{2}\ddc (-\log \omega_D^n/\omega^n)+ \ric \omega= - \theta- \ddc u + [D]=- \ddc u+ \ddc \log |s|_h - \ddc F+ \ric \omega.$$
	Equivalently, we obtain
	\begin{align}\label{eq-KE0Gue}
	\omega_D^n= e^{2u+2F+2c} |s|_h^{-2} \omega^n,
	\end{align}
	for some constant $c>0$. Considering $u+c$ in place of $u$ gives
	\begin{align}\label{eq-ptdinhlky14}
	(\ddc u+ \theta)^n= e^{2u+2F} |s|_h^{-2} \omega^n.
	\end{align}
	By~\cite[Theorem. 4.2]{Berman-Guenancia}, this equation admits a unique solution $u \in \mathcal{E}(X, \theta)$. 
	
	Now consider a constant $\epsilon >0$ small enough so that the class $K_X+ (1-\epsilon) D$ is still big. Observe that 
	$\theta_\epsilon:= \theta- \epsilon \eta$ is  a smooth representative of the class $c_1(K_X+(1-\epsilon)D)$, 
	It was shown in \cite{BEGZ} that there exists a unique K\"ahler--Einstein metric $\omega_\epsilon$ of full Monge-Amp\`ere mass such that 
	\begin{equation}\label{eq: KEmetrics}
	\ric (\omega_\epsilon)= - \omega_\epsilon+ (1-\epsilon)[D], 
	\end{equation}
	or equivalently $\omega_\epsilon= \ddc u_\epsilon+\theta_\epsilon$, where $u_\epsilon \in \mathcal{E}(X, \theta_\epsilon)$ is the unique solution of the equation 
	\begin{align}\label{eq-KE0Gue2}
	(\ddc u_\epsilon+ \theta_\epsilon)^n = e^{2u_\epsilon+ 2F} |s|_h^{-2(1-\epsilon)} \omega^n= e^{2u_\epsilon+ 2F} |s|_h^{-2}(|s|_h^{2\epsilon} \omega^n).
	\end{align}
	Note that although $u_\epsilon - V_{\theta_\epsilon}$ is bounded for every $\epsilon>0$ (see~\cite[Theorem 4.1]{BEGZ}), the function $u- V_\theta$ could be unbounded in general, the structure of $\omega_D$ near $D$ is mostly like a Poincar\'e metrics or a conic one; cf.~\cite{Guenancia-Paun}.
	One sees, hence, that Theorem \ref{the-stability} is more general than Theorem \ref{the-KE-Gue}. We now proceed with the proof of Theorem \ref{the-stability}.

	\begin{remark} In the setting of Theorems~\ref{the-KE-Gue}, by applying Theorem \ref{the-stability} to the equation  (\ref{eq-KE0Gue2}), we see that (\ref{eq-ptdinhlky14}) admits a solution in  $\mathcal{E}(X, \theta)$. This gives an alternative proof of~\cite[Theorem. 4.2]{Berman-Guenancia} without using the variational method.
	\end{remark}

	\begin{remark} Let $(u_\epsilon)_\epsilon$ be the sequence of potentials of the twisted K\"ahler--Einstein metrics $\omega_\epsilon$ in the proof of Theorem \ref{the-TianYaumetric}. We already know that $u_\epsilon$ converges to $u$ (the potential of $\omega_D$) locally in $\mathcal{C}^\infty$-topology in $X \backslash (D \cup E)$. Direct application of Theorem \ref{the-stability} to equations defining $\omega_\epsilon$ implies that $u_\epsilon$ decreases in capacity to $u$. This global property is much stronger than $L^1$-convergence of $u_\epsilon$ to $u$.  
	\end{remark}
	
	\subsection{Discussions} \label{subsec-explain-strategy}
	
	Our goal here is to give an informal discussion about the difficulties in the proof of Theorem \ref{the-stability} (hence of Theorem~\ref{the-KE-Gue}) and explain our strategy to overcome these issues.

	Set $\mu_j:= f_j |s|_h^{-2} \omega^n$. The difficulty in proving the convergence in Theorem \ref{the-stability} is that, unlike some standard situations, the sequence of measures $(\mu_j)_j$ is not convergent in the mass norm (or even in the weak sense). This is due to the fact that $f|s|_h^{-2} \omega^n$ could be of infinite mass on $X$. To be more precise, a usual way to prove the $L^1$ convergence of potentials of $(u_j)_j$ is as follows (see, e.g., \cite{BEGZ}). Let $\delta>0$ be a fixed positive constant. Without loss of generality, we can assume that $u_j$ converges to some $u'$ in $L^1$ (ignore for the moment that $u_j$ might a priori converge to $-\infty$). Hence $u_j$ is $(\theta+ \delta \omega)$-psh for every $j$ big enough. We consider $u'_j:= (\sup_{j' \ge j} u_{j})^*$ which decreases to $u'$. We have  
	$$(\ddc u'_j + \theta+\delta \omega)^n  \ge e^{ \inf_{j' \ge j} u_j} \mu_j.$$ 
	It follows that for every (continuous) function $g \ge 0$, there holds
	$$\liminf_{j \to \infty} \int_X g (\ddc u'_j + \theta+\delta \omega)^n  \ge \liminf_{j \to \infty}\int_X g  e^{ \inf_{j' \ge j} u_j} \mu_\epsilon.$$
	If one had that $\mu_j$ converges to $\mu: =|s|_h^{-2} f\omega^n$ in the mass norm, then it would follow from the last estimate that 
	$$\liminf_{j\to \infty} \int_X g (\ddc u'_j + \theta+\delta \omega)^n  \ge \liminf_{j \to \infty}\int_X g  e^{u'} \mu,$$
	and hence, by letting $\delta \to 0$, 
	$$(\ddc u'+ \theta)^n= e^{u'}\mu,$$
	from this, one would conclude that $u'=u$. However, as already mentioned, it is not true that $\mu_j \to \mu$ in the mass norm because $\mu(X)$ might be equal to infinity.  
	
	Our strategy is to use the quantitative domination principle from \cite{Do-Vu_quantitative} because the standard domination principle is not sufficient for our purpose. Let us explain why.  
	Since $|s|_h^2 \le 1$, we get
	\begin{align}\label{ine-sosanhdominaepsilon}
	( \ddc u_{j_1}+\theta_{j_1})^n \le e^{u_{j_1}-u_{j_2}} ( \ddc u_{j_2}+\theta_{j_2})^n 
	\end{align}
	if $j_1 \le j_2$. So, it is tempting to make use of the domination principles. In order to do so, one needs to consider $u_{j_k}$, $k=1,2$, in the same cohomology class. 
	
	Fix now a constant $\epsilon>0$. Put $\theta'_\epsilon:= \theta+\epsilon \omega$.  Let $j_\epsilon \in \N$ be such that $\theta- \epsilon \omega \le \theta_j \le\theta'_\epsilon$ for $j \ge j_\epsilon$.  Consider $j_\epsilon<j_1 \le j_2$.  Hence  $u_{j_s}$ is $\theta'_\epsilon$-psh functions for $s=1,2$.
	
	\begin{lemma} \label{le-chantrenuKE} We have
		$$-C \epsilon \le \int_X (\ddc V_{\theta_j}+ \theta'_\epsilon)^n- \int_X (\ddc V_{\theta'_\epsilon}+\theta'_\epsilon)^n \le 0,$$
		for some constant $C>0$ independent of $j$.
	\end{lemma}
	
	\proof
	Since $V_{\theta_j} \le V_{\theta'_\epsilon} $, we get 
	$$\int_X (\ddc V_{\theta_j}+ \theta'_\epsilon)^n- \int_X (\ddc V_{\theta'_\epsilon}+\theta'_\epsilon)^n \le 0$$
	by monotonicity of non-pluripolar products. The other desired inequality follows from the fact that $\vol(\{\theta_j\}) \ge \vol(\{\theta- \epsilon \omega\}) \ge \vol({\theta+\epsilon \omega})- O(\epsilon)$. 
	\endproof
	
	By Lemma \ref{le-upperboundforuj}, we can assume that $u_j \le 0$ for every $j\in\N$.
	There are two issues when working with $\theta'_\epsilon$. The first problem is that $u_{j_s}$ is no longer in $\mathcal{E}(X, \theta'_\epsilon)$ for $s=1,2$, but belongs to $\mathcal{E}(X, \theta'_\epsilon, V_{\theta_{j_s}})$, and $ V_{\theta_{j_s}} \le V_{\theta'_\epsilon}$ in general. The second problem is that (\ref{ine-sosanhdominaepsilon}) is no longer true if $\theta_{j_s}$ is replaced by $\theta'_\epsilon$. For a constant $M>0$, direct computations show that on $\{u_{j_1} \le u_{j_2}- M\}$ there holds 
	\begin{align*}
	( \ddc u_{j_1}+\theta'_\epsilon)^n  &\le e^{u_{j_1}-u_{j_2}} ( \ddc u_{j_2}+\theta_{j_2})^n+ O(\epsilon) \\
	& \le e^{-M} ( \ddc u_{j_2}+\theta'_{\epsilon})^n+ O(e^{-M}\epsilon)+ O(\epsilon),
	\end{align*}
	where for a constant $A>0$ we denote by $O(A)$ a measure of mass bounded by $A$. This estimate suggests that the quantitative domination principle can be applied. To this end, one must, however, have that $u_{j_s} \in \mathcal{E}(X, \theta'_\epsilon)$, which is not true in general. For that reason, we consider a big constant $k>0$ and put 
	$$u_{j,k}:= \max\{u_{j}, V_{\theta'_\epsilon}- k\} \in \mathcal{E}(X, \theta'_\epsilon).$$

	\begin{lemma}\label{le-sossanhhaiMA} There is a constant $A>0$ independent of $\epsilon,k$ such that for $j_1,j_2$ big enough, we have
		$$\mathbf{1}_{\{u_{j_1,k} \le u_{j_2,k}- M\}}(\ddc u_{j_1,k}+\theta'_\epsilon)^n \le \mathbf{1}_{\{u_{j_1,k} \le u_{j_2,k}- M\}}e^{-M}(\ddc u_{j_2,k}+ \theta'_\epsilon)^n+ R_{j_1},$$
		where $R_{j_1}$ is a positive measure whose mass is given by
		$$c(j_1,k, \epsilon):= \|R_{j_1}\|=A\epsilon+\int_{\{u_{j_1} \le V_{\theta'_\epsilon} -k \}} (\ddc u_{j_{1},k}+ \theta'_\epsilon)^n.$$
	\end{lemma}

	\proof Observe
	$$\{u_{j_1,k} \le u_{j_2,k}- M\} \subset K:= \{u_{j_1} \le u_{j_2}- M\} \cap \{u_{j_2} > V_{\theta'_\epsilon} -k\}.$$
	Thus
	\begin{align*} 
	\mathbf{1}_K(\ddc u_{j_1,k}+\theta'_\epsilon)^n  & \le \mathbf{1}_{K \cap \{u_{j_1} > V_{\theta'_\epsilon} -k\}} (\ddc u_{j_1}+\theta'_\epsilon)^n+ \mathbf{1}_{K\cap \{u_{j_1} \le V_{\theta'_\epsilon} -k\}} ( \ddc u_{j_1,k}+\theta'_\epsilon)^n\\
	\nonumber
	& \le \mathbf{1}_{K \cap \{u_{j_1} > V_{\theta'_\epsilon} -k\}}  e^{-M} (\ddc u_{j_2,k}+\theta'_\epsilon)^n\\
	&\quad+O(\epsilon)+ \mathbf{1}_{K\cap \{u_{j_1} \le V_{\theta'_\epsilon} -k\}} ( \ddc u_{j_1,k}+\theta'_\epsilon)^n.
	\end{align*}
	This finishes the proof.
	\endproof

	If we want to apply the quantitative domination principle to $u_{\epsilon_j,k}$, we need to check that $c(j_1,k, \epsilon)$ is ``small'', precisely, 
	$$\limsup_{k\to \infty}\limsup_{\epsilon \to 0} \limsup_{j_1 \to \infty}c(j_1,k, \epsilon) =0.$$
Verifying this limit is one of the main difficulties in our proof. We will go into detail in the next sections. 

	\subsection{Low energy estimate}

	We start with a fact about energy.
	The following monotonicity of energy is a direct consequence of \cite[Lemma 3.2]{Vu_Do-MA} (see also \cite{GZ-weighted}).
	
	\begin{lemma} \label{le-mono-energy} Let $\eta$ be a smooth closed $(1,1)$-form in a big cohomology class. Let $u,v \in \mathcal{E}(X,\eta)$ such that $u \le v$. Let $\chi: \R_{\le 0} \to \R_{ \le 0}$ be a convex increasing function such that $\chi(0)=0$. Then we have 
		$$-\int_X \chi(v- V_\eta) (\ddc v+ \eta)^k \wedge \omega^{n-k} \le -2^k \int_X \chi(u- V_\eta) (\ddc u+ \eta)^k \wedge \omega^{n-k}.$$    
	\end{lemma}

	Let $\epsilon>0$ be a small constant. Let $\theta'_\epsilon:= \theta+ \epsilon \omega$. Let $j_\epsilon \in \N$ be such that $\theta_j \le \theta'_\epsilon$ for every $j \ge j_\epsilon$.  Hence $u_j$ is $\theta_\epsilon'$-psh for $j \ge  j_{\epsilon}$. Note that $u_j \in \mathcal{E}(X, \theta_j)$.  Let 
	$$u''_{j,\epsilon}:= \lim_{k\to \infty} P_{\theta'_\epsilon}( \min \{u_j, \ldots, u_k\}),$$
	for $j \ge j_\epsilon$. 
	By Corollary \ref{cor-envelopchanduoi}, we have that $u''_{j,\epsilon}$ is a well-defined $\theta'_\epsilon$-psh function and increases to some $\theta'_\epsilon$-psh function $u''_\epsilon$ as $j\to \infty$ satisfying that  $u''_\epsilon$ decreases to some $\theta$-psh function $u'' \in \mathcal{E}(X, \theta)$ as $\epsilon$ decreases to $0$. 
	
	Let $\tau: \R_{\le 0} \to \R_{\le 0}$ be an increasing convex function such that $\tau(0)=0$ and $\tau(-\infty)= -\infty$ and 
	\begin{align}\label{ine-nangluongcuatauphay}
	-\int_X \tau(u''- V_\theta) (\ddc u''+ \theta)^n < \infty.
	\end{align}
	Observe that $u_j \ge u''_{j,\epsilon}$ for every $j \ge j_\epsilon$. 
	Let 
	$$u_{j,k}:=\max\{u_j,V_{\theta_\epsilon'}-k\}, \quad u''_{j,\epsilon,k}:= \max\{u''_{j,\epsilon,k},V_{\theta_\epsilon'}-k\}$$
	for $k>0$. We define $u''_{\epsilon,k}$ and $u''_k$ respectively for $u''_\epsilon$ and $u''$ similarly.

	\begin{lemma}\label{le0boundedenergy2} We have
		$$\limsup_{j\to \infty}\int_X - \tau(u_{j,k}- V_{\theta_\epsilon'}) (\ddc u_{j,k}+ \theta_\epsilon')^n \le 4^n\int_X - \tau(u''_k- V_{\theta'_\epsilon}) (\ddc u''_k+ \theta'_\epsilon)^n.$$
		for every  $k \ge 0$. 
	\end{lemma}

	\proof We have $u_{j,k} \ge u''_{j,\epsilon,k}$. Hence using 
	Lemma \ref{le-mono-energy} gives
	$$\int_X - \tau(u_{j,k}- V_{\theta_\epsilon'}) (\ddc u_{j,k}+ \theta_\epsilon')^n  \le 2^n \int_X - \tau(u''_{j,\epsilon,k}- V_{\theta_\epsilon'}) (\ddc u''_{j,\epsilon,k}+ \theta_\epsilon')^n.$$
	Letting $j \to \infty$ and using the fact that $ u''_{j,\epsilon,k}$ increases to $u''_{\epsilon,k}$ (hence we get correspondingly the continuity of energy), one gets
	$$\limsup_{j\to \infty} \int_X - \tau(u_{j,k}- V_{\theta_\epsilon'}) (\ddc u_{j,k}+ \theta_\epsilon')^n  \le 2^n \limsup_{j\to \infty}\int_X - \tau(u''_{j,\epsilon,k}- V_{\theta_\epsilon'}) (\ddc u''_{j,\epsilon,k}+ \theta_\epsilon')^n$$
	which is equal to $ 2^n I_{\epsilon,k}$, where
	$$I_{\epsilon,k}:= \int_X - \tau(u''_{\epsilon,k}- V_{\theta_\epsilon'}) (\ddc u''_{\epsilon,k}+ \theta_\epsilon')^n.$$
	Now observe that $u''_{\epsilon,k}$ decreases to $u''_k$. It follows that
	$$I_{\epsilon,k} \le 2^n \int_X - \tau(u''_{k}- V_{\theta_\epsilon'}) (\ddc u''_{k}+ \theta_\epsilon')^n.$$
	The proof is thus complete.
	\endproof

	\begin{lemma}\label{le0boundedenergy2them} 
		There exists a constant $C>0$ such that  we have 
		$$\limsup_{\epsilon \to 0} \int_X -\tau(u''_k- V_{\theta'_\epsilon}) (\ddc u''_k+ \theta'_\epsilon)^n \le C,$$
		for every  $k \ge 0$. 
	\end{lemma}
	

	\proof
	Let $v''_k:= \max\{u'', V_\theta-k\} \ge u''$. We note that $u''_k$ depends on $\epsilon$ but $v''_k$ does not. Furthermore, we have $u''_k= v''_k$ on $\{u''> V_{\theta'_\epsilon}-k\}$ (because $V_\theta \le V_{\theta'_\epsilon}$). Consequently, one can
	decompose
	\begin{align*}
	\int_X - \tau(u''_k- V_{\theta'_\epsilon}) (\ddc u''_k+ \theta'_\epsilon)^n &= \int_{\{u''> V_{\theta'_\epsilon}-k\}} - \tau(u''_k- V_{\theta'_\epsilon}) (\ddc u''_k+ \theta'_\epsilon)^n+\\
	&\quad  \int_{\{u'' \le  V_{\theta'_\epsilon}-k\}} - \tau(u''_k- V_{\theta'_\epsilon}) (\ddc u''_k+ \theta'_\epsilon)^n\\
	&= \int_{\{u''> V_{\theta'_\epsilon}-k\}} - \tau(v''_k- V_{\theta'_\epsilon}) (\ddc v''_k+ \theta'_\epsilon)^n\\
	&\quad  - \tau(-k) \int_{\{u'' \le  V_{\theta'_\epsilon}-k\}}  (\ddc u''_k+ \theta'_\epsilon)^n.
	\end{align*}
	Denote by $I_1, I_2$ the first and the second term on the right-hand side of the last equality. 
	
Recall that $T= \ddc \rho+ \theta$ is a K\"ahler current with analytic singularities in the class of $\theta$, where $\rho \le 0$.	Since $\ddc\rho+\theta\ge \delta\omega$ for $\delta>0$ we see that $\epsilon^{1/2}\rho+(1-\epsilon^{1/2})V_{\theta_\epsilon'}$ is a negative $\theta$-psh function for $\epsilon$ small enough. Since the latter is less than $V_\theta$, we infer 
	\[V_{\theta_\epsilon'}\le (1-\epsilon^{1/2})^{-1}V_\theta-\epsilon^{1/2}(1-\epsilon^{1/2})^{-1}\rho\le V_\theta-2\epsilon^{1/2}\rho,\]
	for $\epsilon>0$ small enough.
	Since $\tau$ is convex and $\tau(0)=0$, one sees that 
	$$\tau(b+c) \ge \tau(b)+ \tau(c)$$
	for every $ b, c\in \R_{\le 0}$. Applying the last inequality to $b:= v''_k- V_{\theta} $ and $c:= V_{\theta}- V_{\theta'_\epsilon}$, one obtains
	$$-\tau(v''_k- V_{\theta'_\epsilon}) \le -\tau(v''_k - V_\theta)- \tau(V_\theta- V_{\theta'_\epsilon}) \le -\tau(v''_k - V_\theta)- \tau(2 \epsilon^{1/2} \rho).$$
	It follows that 
	\begin{align*}
	I_1 &\le \int_X - \tau(v''_k- V_{\theta}) (\ddc v''_k+ \theta'_\epsilon)^n - \int_X \tau(2 \epsilon^{1/2} \rho)(\ddc u''_k+ \theta'_\epsilon)^n
	\end{align*}
	Letting $\epsilon \to 0$ gives
	\begin{align}\label{ine-danhgiaI1}
	\limsup_{\epsilon \to 0}I_1 \le \int_X - \tau(v''_k- V_{\theta}) (\ddc v''_k+ \theta)^n \le 2^n \int_X - \tau(u''- V_{\theta}) (\ddc u''+ \theta)^n 
	\end{align}
	because $\tau(0)=0$ and $(\ddc u''_k+ \theta)^n$ is a Monge--Amp\`ere measure of bounded potentials. 
	
	We treat $I_2$. Direct computations show
	\begin{align}\label{ine-I2nhe}
	I_2 &= -\tau(-k) \int_X (\ddc u''_k + \theta'_\epsilon)^n + \tau(-k) \int_{\{u''> V_{\theta'_\epsilon}-k\}} (\ddc u''_k + \theta'_\epsilon)^n\\
	\nonumber
	&=  -\tau(-k) \int_X (\ddc V_{\theta'_\epsilon} + \theta'_\epsilon)^n + \tau(-k) \int_{\{u''> V_{\theta'_\epsilon}-k\}} (\ddc u'' + \theta'_\epsilon)^n\\
	\nonumber
	&\le  -\tau(-k) \bigg( \vol(\{\theta'_\epsilon\})- \vol(\{\theta\})+ \vol(\{\theta\})- \int_{\{u''> V_{\theta'_\epsilon}-k\}} (\ddc u'' + \theta)^n\bigg)\\
	\nonumber
	&= -\tau(-k) \bigg( \vol(\{\theta'_\epsilon\})- \vol(\{\theta\})+  \int_{\{u'' \le  V_{\theta'_\epsilon}-k\}} (\ddc u'' + \theta)^n\bigg)\\
	\nonumber
	& \le -\tau(-k) \bigg( \vol(\{\theta'_\epsilon\})- \vol(\{\theta\})+  |\tau(-k)|^{-1} \int_{\{u'' \le  V_{\theta'_\epsilon}-k\}} -\tau(u''- V_{\theta'_\epsilon}) (\ddc u'' + \theta)^n\bigg)\\
    & \le -\tau(-k) \big( \vol(\{\theta'_\epsilon\})- \vol(\{\theta\}) \big)+\int_X-\tau(u''- V_{\theta'_\epsilon}) (\ddc u'' + \theta)^n \nonumber
	\end{align}
	whose limsup as $\epsilon \to 0$ is bounded by a constant independent of $k$, due to~\eqref{ine-nangluongcuatauphay}. 
	Combining (\ref{ine-I2nhe}) and (\ref{ine-danhgiaI1}) gives the desired assertion. 
	\endproof
	
	Combining Lemmas \ref{le0boundedenergy2them} and \ref{le0boundedenergy2}, we obtain the following crucial estimate.

	\begin{proposition}\label{pro-ernergyestimate} 
		There exists a constant $C>0$ independent of $k$ such that  we have
		$$\limsup_{\epsilon \to 0}\limsup_{j\to \infty}\int_X -\tau(u_{j,k}- V_{\theta_\epsilon'}) (\ddc u_{j,k}+ \theta_\epsilon')^n \le C.$$
	\end{proposition}

	\begin{corollary}\label{cor-chantrenphandudomination}
		Let $c(j_1,k,\epsilon)$ be the term defined in Lemma \ref{le-sossanhhaiMA}. Then we have 
		$$\limsup_{k\to \infty}\limsup_{\epsilon \to 0} \limsup_{j_1 \to \infty}c(j_1,k, \epsilon) =0$$
	\end{corollary}
	
	\proof
	Observe
	\begin{align*}
	\limsup_{\epsilon \to 0} \limsup_{j_1 \to \infty}c(j_1,k, \epsilon) &= \limsup_{\epsilon \to 0} \limsup_{j_1 \to \infty} \int_{\{u_{j_1,k} \le V_{\theta'_\epsilon}-k\}} (\ddc u_{j_1,k}+ \theta'_\epsilon)^n\\
	& \le  |\tau(-k)|^{-1} \limsup_{\epsilon \to 0} \limsup_{j_1 \to \infty} \int_{\{u_{j_1} \le V_{\theta'_\epsilon}-k\}}-\tau(u_{j_1,k}- V_{\theta'_\epsilon}) (\ddc u_{j_1,k}+ \theta'_\epsilon)^n\\
	& \le |\tau(-k)|^{-1}C
	\end{align*}
	by Proposition \ref{pro-ernergyestimate}. Letting $k \to\infty$ gives the desired equality. 
	\endproof

	We now finish the proof of Theorem \ref{the-stability}. 
	
	\begin{proof}[End of the proof of Theorem \ref{the-stability}] 
		We apply the quantitative domination principle to $u_{j, k}$ in the class of $\{\theta'_\epsilon\}$ and $\tilde{\chi}:= \tau$ (we can certainly choose $\tau$ so that $\tau(-1)=-1$).
		Let 
		$$B(j_1,j_2,k,\epsilon):= -\int_X \tau(u_{j_1,k}- V_{\theta'_\epsilon}) (\ddc u_{j_1,k}+ \theta'_\epsilon)-\int_X \tau(u_{j_2,k}- V_{\theta'_\epsilon}) (\ddc u_{j_2,k}+ \theta'_\epsilon).$$
		By Proposition~\ref{pro-ernergyestimate}, one has
		\begin{align} \label{ine-chanB}
		\limsup_{\epsilon \to 0}\limsup_{j_1,j_2\to \infty} B(j_1,j_2,k,\epsilon) \le C,
		\end{align}
		where $C>0$ is a constant independent of $k$. Let $0<\kappa \le 1/4$ be a constant. By Theorem~\ref{the domination} applied to $u_{j_1,k}, u_{j_2,k}- \kappa$ (noticing also Lemma \ref{le-sossanhhaiMA}), there exists a constant $C_1>0$ independent of $j_1,j_2,k,\epsilon,\kappa$ such that  
		$$\capK_\omega \big(u_{j_1, k} - u_{j_2, k} \le -2\kappa \big) \le  C_1\kappa^{-2}\frac{\big(B(j_1,j_2,k,\epsilon)\big)^2+1}{h^{\circ n}\big(1/c(j_1,k,\epsilon)\big)},$$
		with $h(t):= (-\tau(-t))^{1/2}$, and $j_1,j_2>j_{\epsilon}$ big enough. Since 
		$$\capK_\omega(u_{j}< V_{\theta'_\epsilon}-k) \le \capK_\omega(u_{j}< -k) \lesssim k^{-1},$$
		one infers
		$$\capK_\omega \big(u_{j_1} - u_{j_2} \le -2\kappa\big) \lesssim \kappa^{-2} \frac{\big(B(j_1,j_2,k,\epsilon)\big)^2+1}{h^{\circ n}\big(1/c(j_1,k,\epsilon)\big)}+ k^{-1}.$$
		Hence 
		$$\limsup_{j_1,j_2 \to \infty} \capK_\omega \big(u_{j_1} - u_{j_2}\le -\kappa\big) \lesssim \kappa^{-2}\limsup_{j_1,j_2 \to \infty} \frac{\big(B(j_1,j_2,k,\epsilon)\big)^2+1}{h^{\circ n}\big(1/c(j_1,k,\epsilon)\big)}+ k^{-1}$$
		for every $\epsilon \in (0,1]$ and $k \in \N$.  Letting $\epsilon \to 0$ and using (\ref{ine-chanB}), we get 
		$$\limsup_{j_1,j_2 \to \infty} \capK_\omega \big(u_{j_1} - u_{j_2}\le -\kappa\big) \lesssim \kappa^{-2}\frac{1}{h^{\circ n}\big(1/c(k)\big)}+ k^{-1},$$
		for every $k$, 
		where 
		$$c(k):= \limsup_{\epsilon \to 0}\limsup_{j_1\to \infty}c(j_1,k,\epsilon).$$
		By Corollary \ref{cor-chantrenphandudomination}, we get $\limsup_{k\to\infty}c(k)=0$. 
		It follows that 
		$$\limsup_{j_1,j_2 \to \infty} \capK_\omega \big(u_{j_1} - u_{j_2}\le -\kappa\big) \lesssim \limsup_{k\to \infty}\kappa^{-1}\frac{1}{h^{\circ n}\big(1/c(k)\big)}+ k^{-1} = 0,$$
		for every constant $\kappa>0$. In other words, the sequence $(u_j)_j$ is decreasing in capacity. Let $\psi_\gamma$ be the function defined right before Corollary \ref{cor-subsolution}. Since $u_j \ge \psi_\gamma$  which is locally bounded outside the union $W$ of $D$ and the non-K\"ahler locus of $\{\theta\}$, we see that $(u_j)_j$ 
		admits a subsequence converging in $L^1$. Let $u$ be, hence, a $L^1$-limit of $(u_j)_j$. Corollary \ref{cor-subsolution} tells us that $u \in \mathcal{E}(X, \theta)$.


		 Since $u_j \ge \psi_\gamma$, we also get $u \ge \psi_\gamma$ and hence $u$ is locally bounded outside $W$. Proposition \ref{prop: cv_decreasingcapacity} now implies that 
		$$(\ddc u_j + \theta_j)^n \to (\ddc u+ \theta)^n$$
		weakly  on $X \backslash W$ as $j \to \infty$. It follows that 
		$$(\ddc u+ \theta)^n= e^u |s|_h^{-2} f \omega^n$$
		on $X \backslash W$. The equality holds indeed on $X$ because the non-pluripolar products have no mass on $W$, which is a pluripolar set. Hence we get $(\ddc u+ \theta)^n= e^u |s|_h^{-2} f \omega^n$  on $X$ and $u \in \mathcal{E}(X,\theta)$. This is, at most, one such $u$ by the domination principle. This yields that $u_j \to u$ in $L^1$ and $u$ satisfies the required Monge--Amp\`ere equation.  This ends the proof of Theorem~\ref{the-stability}.
	\end{proof}


	\bibliography{biblio_family_MAt,bib-kahlerRicci-flowt}
	\bibliographystyle{siam}
	
	\bigskip

	\noindent\Addresses
\end{document}